\documentclass[12pt]{article}

\usepackage{vmargin}
\usepackage{fancyhdr}
\usepackage{ifthen}
\usepackage{graphicx} 
\usepackage{indentfirst}
\usepackage[english]{babel}

\usepackage{fouriernc}
\DeclareMathAlphabet{\mathcal}{OMS}{cmsy}{m}{n} 

\usepackage{amsmath}   
\usepackage{amssymb}   
\usepackage{amsfonts}
\usepackage{amsthm}
\usepackage{mathrsfs}

\usepackage{pst-node} 
\usepackage{pstricks}

\usepackage{tweaklist} 
\renewcommand{\itemhook}{\setlength{\topsep}{1pt}%
\setlength{\itemsep}{-1pt}}

\renewcommand{\enumhook}{\setlength{\topsep}{1pt}%
\setlength{\itemsep}{-1pt}}

\usepackage[colorlinks=true,linkcolor=magenta,citecolor=blue]{hyperref}

\usepackage{comment}
\usepackage{color}

\specialcomment{notes}{%
  \begingroup\color{red}\begin{center}\begin{minipage}{12cm} $\blacktriangleright\ $}{%
  \hfill$\ \blacktriangleleft$ \end{minipage} \end{center}\endgroup}

\setpapersize{A4}
\setmarginsrb{33mm}{20mm}{33mm}{25mm}{10mm}{6mm}{0mm}{10mm}
\fancypagestyle{plain}{%
\fancyhf{} 
\lfoot[\thepage]{}
\cfoot[]{}
\rfoot[]{\thepage}
}

\newcommand{\sk}{\smallskip}
\newcommand{\mk}{\medskip}
\newcommand{\bk}{\bigskip}

\newcommand{\G}{\mathbf{G}}
\newcommand{\B}{\mathbf{B}}
\newcommand{\U}{\mathbf{U}}

\newcommand{\T}{\mathbf{T}}
\renewcommand{\P}{\mathbf{P}}
\renewcommand{\L}{\mathbf{L}}

\newlength{\leftlength}
\newlength{\rightlength}
\newlength{\calculskip}
\newcommand{\calculvskip}[1]{%
  \ifthenelse{#1 = 0}{\setlength{\calculskip}{0pt}}{}%
  \ifthenelse{#1 = 1}{\setlength{\calculskip}{\smallskipamount}}{}%
  \ifthenelse{#1 = 2}{\setlength{\calculskip}{\medskipamount}}{}%
  \ifthenelse{#1 = 3}{\setlength{\calculskip}{\bigskipamount}}{}%
  \ifthenelse{#1 = 4}{\setlength{\calculskip}{1cm}}{}%
  \vskip\calculskip
}

\newcommand{\leftcentersright}[4][2]{%
        \settowidth{\leftlength}{#2}%
        \settowidth{\rightlength}{#4}%
        \calculvskip{#1}
        \noindent#2\hskip-\leftlength%
        \hfill#3\hfill
        \mbox{}\hskip-\rightlength#4%
        \vskip\calculskip%
        }

\newcommand{\centers}[2][2]{\leftcentersright[#1]{}{#2}{}}
\newcommand{\leftcenters}[3][2]{\leftcentersright[#1]{#2}{#3}{}}

\makeatletter
\def\svhline{%
  \noalign{\ifnum0=`}\fi\hrule \@height2\arrayrulewidth \futurelet
   \reserved@a\@xhline}
   
\def\hlinewd#1{%
\noalign{\ifnum0=`}\fi\hrule \@height #1 %
\futurelet\reserved@a\@xhline}
 \makeatother

\numberwithin{equation}{section}

\newtheorem{prop}[equation]{Proposition}  

\newtheorem{thm}[equation]{Theorem}

\newtheorem*{thm2}{Theorem}

\theoremstyle{definition}

\newtheorem{rmk}[equation]{Remark}

\newtheorem{hyp}[equation]{Assumption}

\newcommand{\X}{\mathrm{X}}

\newcommand{\Z}{\mathrm{Z}}
\newcommand{\Hc}{\mathrm{H}_c}

\newcommand{\Rgc}{\mathrm{R}\Gamma_c}

\newcommand{\F}{\mathbb{F}}
\newcommand{\qlb}{\overline{\mathbb{Q}}_\ell}
\renewcommand{\H}{\mathrm{H}}
\usepackage{pstricks-add}
\usepackage{lscape}
\usepackage{multirow}
\usepackage{rotating}

\begin{document}

\title{Cohomology of Deligne-Lusztig varieties for short-length regular elements in exceptional groups}
\author{Olivier Dudas\footnote{Oxford Mathematical Institute.}
\footnote{The author is supported by the EPSRC, Project No EP/H026568/1 and by Magdalen College, Oxford.}}

\maketitle

\begin{abstract} We determine the cohomology of Deligne-Lusztig varieties associated to some short-length regular elements for split groups of type $F_4$ and $E_n$. As a byproduct, we obtain conjectural Brauer trees for the principal $\Phi_{14}$-block of $E_7$ and the principal $\Phi_{24}$-block of $E_8$.
\end{abstract}

\section*{Introduction}

Let $G$ be a finite group and $\ell$ be a prime number. The $\ell$-modular representation theory of $G$ is somehow controlled by the representation theory of local subgroups, namely the $\ell$-subgroups of $G$ and their normalisers. Broué's abelian defect conjecture is one of the major open problems in this framework: it predicts that  an $\ell$-block of $G$ with abelian defect group is derived equivalent to its Brauer correspondent. From the work of Rickard \cite{Ri4}, we know that such an equivalence should be induced by a perfect complex. Unfortunately, there is no canonical construction in general. 

\sk

When $G =\G^F$ is a finite reductive group, Broué's  suggested that the complex representing the cohomology of some Deligne-Lusztig variety should be a good candidate. Together  with Michel in \cite{BMi2}, they made precise which specific Deligne-Lusztig varieties would be associated to principal $\Phi_d$-blocks when $d$ is a regular number. They introduced the notion of \emph{good $d$-regular elements} $w \in W$ and conjectured that

\begin{itemize}

\item for $i \neq j$, the groups $\Hc^i(\X(w),\qlb)$ and $\Hc^j(\X(w),\qlb)$ have no irreducible constituents in common;

\item the irreducible constituents of $\Hc^\bullet(\X(w),\qlb)$ are exactly the unipotent characters lying in the principal $\Phi_d$-block;

\item the endomorphism algebra $\mathrm{End}_{\G^F}\big(\Hc^\bullet(\X(w),\qlb)\big)$ is a \emph{$d$-cyclotomic Hecke algebra}.

\end{itemize}

\noindent As for now, these statements have been veryfied in very few cases only. Computing the cohomology of a Deligne-Lusztig variety is a difficult problem, and the only results in this direction have been obtained by Lusztig in \cite{Lu} when $d$ is the Coxeter number (that is when $w$ is a Coxeter element of $W$), for groups of rank $2$ by Digne, Michel and Rouquier in \cite{DMR} and for $d=n$ in type $A_n$ and $d=4$ in type $D_4$ by Digne and Michel in \cite{DM2}. The purpose of this paper is to provide new examples for exceptional groups and in the spirit of Broué's conjecture, to deduce structural properties of the corresponding $\Phi_d$-block.

\sk

We shall adapt Lusztig's strategy: if a character is non-cuspidal then it should appear in the cohomology of a certain quotient of the Deligne-Lusztig variety $\X(w)$. In the Coxeter case, Lusztig proved that this quotient can be expressed in terms of a Deligne-Lusztig variety associated to a "smaller" Coxeter element, providing an inductive method to compute the cohomology of $\X(w)$. The first step towards our main result is to show an analogous property for the $d$-regular elements we are interested in. To this end we will make  extensive use of \cite{Du5}. Unfortunately, this will not give enough information to deal with non-principal series.
In order to compute the cuspidal part of the cohomology of $\X(w)$, we shall, as in \cite{Lu}, introduce compactifications of $\X(w)$. Unlike the Coxeter case, the cuspidal part of $\Hc^\bullet (\X(w),\qlb)$ might not be concentrated in degree $\ell(w)$ since some of the divisors of $\overline{\X}(w)$  might provide cuspidal characters. However, the results in \cite{DMR} are sufficient to determine explicitely in which groups they actually appear and we obtain the following result:

\begin{thm2} Let $w$ be a good $d$-regular element. Then the contribution of the principal series and the discrete series to the cohomology of the Deligne-Lusztig variety $\X(w)$ is explicitely known in the following cases:

\begin{itemize}

\item $(\G,F)$ has type $F_4$ and $d=8$;

\item $(\G,F)$ has type $E_6$ and $d=9$;

\item $(\G,F)$ has type $E_7$ and $d=14$;

\item $(\G,F)$ has type $E_8$ and $d=24$.

\end{itemize}
\end{thm2}

\noindent We will also obtain partial results for the other series, as well as predictions on their contribution, in line with the formula given by Craven in \cite{Cra}.

\sk

Using Lusztig's results in the Coxeter case, Hiss, L\"ubeck and Malle have conjectured that the Brauer tree of the principal $\Phi_h$-block can be read off the cohomology of the Coxeter variety \cite{HLM}. Using the existing Brauer trees and the previous theorem, we propose conjectural planar embedded Brauer trees for the principal $\Phi_{14}$-block of $E_7$ and for the principal $\Phi_{24}$-block of $E_8$ (see Figure \ref{E7} and \ref{E8}). We believe that a further study of the cohomology of the corresponding Deligne-Lusztig varieties as  in \cite{Du3}, \cite{Du4} and \cite{DR} should  give credit to these predictions.

\section{Methods for determining the cohomology}

Let $\G$ be a connected reductive group, together with a Frobenius $F$ defining a  $\F_q$-structure on $\G$. If $\mathbf{H}$ is any $F$-stable algebraic subgroup of $\G$, we will denote by $H$ the finite group of fixed points $\mathbf{H}^F$.  We fix a Borel subgroup $\B$ containing a maximal torus $\T$ of $\G$ such that both $\B$ and $\T$ are $F$-stable. The associated Weyl group is $W = N_\G(\T)/\T$ and the set of simple reflections will be denoted by $S$. We will assume that $(\G,F)$ is split, so that $F$ acts trivially on $W$. 

\sk 
Recall from \cite{DeLu} that to any element $w\in W$ one can associate a \emph{Deligne-Lusztig variety} 

\centers{$\X(w) \, = \, \big\{ g \B \in \G/\B \, | \, g^{-1} \, {}^F g \in \B w \B\big\}.$}

\noindent It is a quasi-projective variety of pure dimension $\ell(w)$, on which $G$ acts by left multiplication. This definition has been subsequently generalized in \cite{BMi2} to elements of the Artin-Tits monoid $B^+$.

\sk

The $\ell$-adic cohomology of $\X(w)$ carries a lot of information on  ordinary and modular representations of $G$. Throughout this paper, we will be interested in the case where $w$ is a good $d$-regular element, or equivalently when $w$ is a $d$-root of $\boldsymbol \pi = \mathbf{w}_0^2$ in the Braid group of $W$. In that case, it is conjectured that the cohomology of $\X(w)$ gives a good model for the unipotent part of the principal $\Phi_d$-block (see for example \cite{BMi2} and \cite{BMM} or the introduction for more details).

\subsection{Non-cuspidal characters}

To any subset $I \subset S$ one can associate a standard parabolic subgroup $\P_I$ containing $\B$ and a standard Levi subgroup $\L_I$ containing $\T$. If $\U_I$ denotes the unipotent radical of $\P_I$, the parabolic subgroup decomposes as $\P_I = \L_I \U_I$ and both $\L_I$ and $\U_I$ are $F$-stable. By \cite[XVII, 6.2.5]{SGA4}, the $U_I$-invariant part of the cohomology of $\X(w)$ is isomorphic to the cohomology of $U_I \backslash \X(w)$. Consequently, one can detect the presence of non-cuspidal modules in the cohomology of $\X(w)$ by studying the quotient variety  $U_I \backslash \X(w)$ for various subsets $I$. In some specific cases, we can express such a quotient (or at least its cohomology) by means of smaller Deligne-Lusztig varieties. 

\sk 

Let  $\mathbf{b} = \mathbf{w}_1 \cdots \mathbf{w}_r \in B^+$ be an element of the Braid monoid, written as a product of reduced elements (i.e. $w_i \in W$). Recall from \cite{Du5} that the decomposition of $\G/\B$ into $\P_I$-orbits induces a decomposition of $\X(\mathbf{b})$ into locally closed $P_I$-subvarieties, called \emph{pieces}

\centers{$\X_{(W_Ix_1,\ldots,W_Ix_r)}(\mathbf{b}) \, = \, \X(\mathbf{b}) \cap \big(\P_I x_1 \B/\B \times \dots \times \P_I x_r \B/\B\big)$}

\noindent where each $x_i$ runs over the set of $I$-reduced elements of $W$. When $I$ and $\mathbf{b}$ are clear from context, we shall simply denote this variety by $\X_{(x_1,\ldots, x_r)}$. Throughout this paper, we will make extensive use of a particular case of the main theorem of \cite{Du5}. It can be deduced from \cite[Remark 3.12]{Du5} when each set $I_i$ is empty and all the $x_i$'s are equal to the same element:

\begin{thm}\label{mainthm}Let $\mathbf{b} = \mathbf{w}_1 \cdots \mathbf{w}_r \in B^+$ with $w_i \in W$, $I$ be a subset of $S$ and $x$ be an $I$-reduced element of $W$. We assume that each  $w_i$ can be decomposed as $w_i = \gamma_i w_i'$ with $\gamma_i \in S \cup \{1\}$ and $w_i' \leq w_i$ be such that

\begin{itemize}

\item[$\mathrm{(a)}$] if $\gamma_i = 1$ then $v_i = x w_i x^{-1} \in W_I$;

\item[$\mathrm{(b)}$] if $\gamma_i \in S$ then $x\gamma_i x^{-1} \notin W_I$, $v_i = x w_i' x^{-1} \in W_I$ and $\ell(w_i') = \ell(v_i)$.

\end{itemize}

\noindent Let $d$ be the number of $w_i$'s satisfying condition $\mathrm{(b)}$ and $ e = \sum \dim (\U_I^x \cap {}^{w_i'} \U \cap \U^-)$. Then we have the following isomorphism of graded $L_I \times \langle F \rangle$-modules:

\centers{$ \Hc^\bullet(U_I \backslash \X_{(x,\ldots, x)}) \, \simeq \, \Hc^\bullet\big((\mathbb{G}_m)^d \times \X_{\L_I} (\mathbf{v}_1 \cdots \mathbf{v}_r)\big)[-2e](e).$}

\end{thm}

\begin{rmk}\label{rmk1}In the particular cases we will be interested in, $\mathbf{b}$ will always be reduced. In that case, it corresponds to an element $w \in W$ and we have $w = w_1 \cdots w_r$ with $\ell(w) = \ell(w_1) + \cdots + \ell(w_r)$. Note that in general the variety $\X(\mathbf{b}) \subset (\G/\B)^r$ can have much more pieces that $\X(w) \subset \G/\B$, since

\centers{$ \X_{W_Ix}(w) \, = \, \displaystyle \bigcup_{\begin{subarray}{c} x_2, \ldots,x_r \\ I\text{-reduced}\end{subarray}} \X_{(W_I x , W_Ix_2,\ldots , W_I x_r)}(\mathbf{b}).$}

\noindent However, in our specific examples we will observe that the piece $\X_{(W_I x , W_Ix_2,\ldots , W_I x_r)}$ will be empty unless $x_2 = \cdots = x_r = x$, so that $\X_x \simeq \X_{(x,x,\ldots, x)}$.
\end{rmk}

\subsection{Cuspidal characters\label{5se22}}

By definition, cuspidal representations of $G$ have no non-zero element invariant under the action of $U_I$ unless $I = S$. In particular, the cohomology of the quotient variety $U_I \backslash \X(w)$ contains no information on the cuspidal characters that can appear in $\X(w)$. In this section we shall briefly review some methods developped in \cite{DM2} and \cite{DMR} in order to solve the problem of finding cuspidal characters in the cohomology of Deligne-Lusztig varieties.

\sk

Let $\mathbf{b} = \mathbf{w}_1 \cdots \mathbf{w}_r$ with $w_i \in W$. Recall that the variety $\X(\mathbf{b})$ has a nice compactification 

\centers{$  \X(\underline{w}_1 \cdots \underline{w}_r) =  \big\{ (g_0,g_1,\ldots, g_r) \in (\G/\B)^{r+1} \, \big| \, g_{i-1}^{-1} g_{i} \in \overline{\B w_i \B}\  \text{and}  \ g_{r}^{-1} F(g_{0}) \in \B \big\}$}

\noindent  which has the following properties (see \cite{DM2} and \cite{DMR}) :

\begin{prop}\label{5prop3}Let $w_1, \ldots, w_r$ be elements of $W$, 

\begin{enumerate} 

\item[$\mathrm{(i)}$] $\X(\underline{w}_1 \cdots \underline{w}_r)$ is a projective variety of dimension  $\ell(w_1) + \cdots + \ell(w_r)$;

\item[$\mathrm{(ii)}$] $\X(\underline{w}_1 \cdots \underline{w}_r )$ is smooth whenever each variety  $\overline{\B w_i \B}$ is;

\item[$\mathrm{(iii)}$] $\X(\underline{w}_1 \cdots \underline{w}_r)$ is rationally smooth whenever each variety $\overline{\B w_i \B}$ is;

\item[$\mathrm{(iv)}$] $\X(\underline{w}_1 \cdots \underline{w}_r)$ has a filtration by closed subvarieties $\X(\underline{v}_1 \cdots \underline{v}_r)$ where the $v_i$'s satifisfy $v_i \leq w_i$.

\end{enumerate}

\end{prop}

\begin{rmk}\label{5rem2}A particular case is when each $w_i$ is a simple reflection $s_i$. Then the variety $\X(\underline{w}_1 \cdots \underline{w}_r)$ coincides with the smooth compactification introduced by Deligne and Lusztig in \cite{DeLu}.
\end{rmk}

Let $w \in W$. In order to compute the cuspidal part of the cohomology of $\X(w)$ using the previous compactifications, we will use the following results:

\begin{itemize}

\item[$\mathbf{(C1)}$] \label{cusp1}the cohomology of $\X(w)$ over $\qlb$ is zero outside the degress $\ell(w), \ldots, 2 \ell(w)$ \cite[Corollary 3.3.22]{DMR};

\item[$\mathbf{(C2)}$]  \label{cusp2}the following triangle is distinguished in $D^b(\qlb G$-$\mathrm{Mod})$ :

\centers{$ \Rgc\big(\X(w), \qlb\big)\ \longrightarrow  \ \Rgc\big(\X(\underline{w}), \qlb\big) \ \longrightarrow \Rgc\Big( \displaystyle \bigcup_{v < w} \X(\underline{v}) , \qlb \Big)  \ \rightsquigarrow $}

\item[$\mathbf{(C3)}$]  \label{cusp3}when $\X(\underline{w})$ is rationally smooth,  its cohomology as a graded $G \times \langle F\rangle_{\mathrm{mon}}$-module  can be explicitely computed using \cite[Corollary 3.3.8]{DMR};

\item[$\mathbf{(C4)}$]  \label{cusp4}let $\rho$ be a cuspidal representation of $G$ that appears in the cohomology of a Deligne-Lusztig variety associated to a Coxeter element of $W$. If $w$ itself is not a Coxeter element, any eigenvalue $\lambda$ of $F$ on $\Hc^{\ell(w)}(\X(w),\qlb)_\rho$ satisfy $|\lambda| <  |q^{\ell(w)/2}|$. This is a particular case of \cite[Proposition 3.3.21]{DMR}. 

\end{itemize}

\sk

Note finally that the property of being rationally smooth of $\X(\underline{w})$ can be read off the Kazdhan-Lusztig polynomials of $W$ \cite[proposition 3.2.5]{DMR}. If $\X(\underline{w})$ happens to be not rationally smooth, we can always decompose $w$ into a product $w=w_1 \cdots w_r$ such that each variety $\overline{\B w_i \B}$ is.

\section{Some particular cases}

For short-length regular elements, one can observe that only a small number of pieces $\X_x$ are non-empty. In addition, they very often satisfy the assumptions of Theorem \ref{mainthm}. For some of these elements, we can therefore compute explicitely the cohomology of the quotient  $U_I \backslash \X(w)$, and eventually deduce the cohomology of $\X(w)$ using the results discussed in Section \ref{5se22}.

\sk

To make the computations easier, we shall use the notation introduced in \cite{DM2}: the cohomology of the Deligne-Lusztig variety $\X(w)$ as a graded $G \times \langle F \rangle_{\mathrm{mon}}$-module will be represented by a polynomial $\H_{\X(w)}(t^{1/2},h)$ with coefficients in the semi-group $\mathbb{N} \, \mathrm{Irr}\, G$.  By a theorem of Lusztig, when $\rho$ is a unipotent character, any eigenvalue of $F$ on the $\rho$-isotypic part of $\Hc^i(\X(w),\qlb)$ can be written $\lambda_\rho q^{ j/2}$, where $\lambda_\rho$ is a root of unity independent of $w$ and $i$. The multiplicity of $\rho$ in  $\Hc^i(\X(w),\qlb)$ with eigenvalue $\lambda_\rho q^{ j/2}$ will be encoded by the coefficient of $h^i t^{j/2}$ in the polynomial $\H_{\X(w)}(t^{1/2},h)$. For example, if $\X(s)$ is the Drinfeld curve for $\G = \mathrm{SL}_2(\overline\F_p)$ then $\H_{\X(w)} = h \mathrm{St} + h^2 t \mathrm{Id}$. 

\sk

Since we are dealing with exceptional Weyl groups, and more specifically with the combinatorics of distinguished subexpressions, we will use the package $\mathsf{CHEVIE}$ in $\mathsf{GAP}$. We have  written a couple of useful functions to determine whether a piece of a Deligne-Lusztig variety is non-empty, and to describe it that case its quotient by a finite unipotent group. These functions can be found in \cite{Gap} (or will be soon available) under the name $\mathsf{deodhar.g}$.

\subsection{\texorpdfstring{$8$}{8}-regular elements for groups of type \texorpdfstring{$F_4$}{F4}}

Let $(\G,F)$ be a split group of type $F_4$. To fix the notation we will consider the following Dynkin diagram:

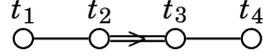
\centers{\begin{pspicture}(3,0.7)

  \cnode(0,0.05){4pt}{A}
  \cnode(1,0.05){4pt}{B}
  \cnode(2,0.05){4pt}{C}
  \cnode(3,0.05){4pt}{D}

  \ncline[nodesep=0pt]{A}{B}\naput[npos=-0.2]{$t_1$}\naput[npos=1.2]{$t_2$}
  \ncline[nodesep=0pt]{C}{D}\naput[npos=-0.2]{$t_3$}\naput[npos=1.2]{$t_4$}
  \ncline[nodesep=0pt,doubleline=true]{B}{C}\ncput{$>$}
\end{pspicture}}

\noindent where $t_1,t_2,t_3$ and $t_4$ are the simple reflections.

\sk

Recall that there exist $d$-regular elements for $d \in\{1,2,3,4,6,8,12\}$ only (see for example \cite{Sp}). Note that the largest integer corresponds to the Coxeter number.  By \cite{BMi2}, for any such $d$ one can find a particular $d$-regular element which is a $d$-th root of $\boldsymbol \pi$ in the Braid monoid. By \cite[11.22]{Be} and \cite[Proposition 5.5]{DM3}, the cohomology of the corresponding Deligne-Lusztig variety does not depend on the choice of a root. For our purposes we will take

\centers{$ w = t_1 t_2 t_3  t_2  t_3 t_4. $}

\mk

\noindent \textbf{\thesubsection.1.  Cohomology of $U_I \backslash \X(w)$.}\label{211} We start by computing the cohomology of the quotient  $U_I \backslash \X(w)$ where $I = \{t_2,t_3\}$. Using the criterion given in  \cite[Lemma 8.3]{DM2} and the package $\mathsf{CHEVIE}$ in $\mathsf{GAP}$ one can check that there are only three non-empty pieces $\X_x$, corresponding  to the cosets $W_I x = W_Iw_0$, $W_Iw_0t_1t_2$ and $W_Iw_0t_4 t_3$.  Theorem \ref{mainthm} does not apply directly to all of these cells, but we can add an intermediate step. Let $J = \{t_2,t_3,t_4\}$ and $K = \{t_1,t_2,t_3\}$. We have $W_J w_0 = W_J w_0 t_4 t_3$ (resp. $W_K w_0 = W_K w_0 t_1 t_2$) whereas $\X_{W_I w_0 t_1 t_2} = \X_{W_J w_0 t_1 t_2}$  (resp.  $ \X_{W_J w_0 t_4 t_3}$) is stable by $P_J$ (resp. $P_K$). Therefore only two pieces appear in the decompositon of $U_J \backslash \X(w)$ (resp. $U_K \backslash \X(w)$). 

\begin{itemize}

\item Let $y$ be the minimal element of $W_J w_0 t_1 t_2$. Since $t_1t_2$ is $J$-reduced, $y = w_J w_0 t_1 t_2$. Let us decompose $w$ as $w = w_1 w_2 $ with $w_1 = t_1 t_2 t_3 t_2$ and $w_2 = t_3 t_4 = t_3 w_2'$. We have ${}^y w_1 = t_2 t_3 \in W_J$ and ${}^y w_2' = t_4$ and therefore we can apply Theorem \ref{mainthm} to compute the cohomology of the piece of $\X(\mathbf{w}_1 \mathbf{w}_2)$ corresponding to $(W_J y, W_Jy)$. Furthermore, one can check (using $\mathsf{GAP}$ again) that the pieces corresponding to $(W_Jy,W_Jy')$ are empty unless $y$ and $y'$ lie in the same coset. In particular, $\X_{W_Jy}(w) \simeq \X_{W_J y, W_Jy} (\mathbf{w}_1 \mathbf{w}_2)$ and

\centers{$  \Hc^\bullet (\X_{W_J y}, \qlb)^{U_J} \, \simeq \, \Hc^\bullet(\mathbb{G}_a \times  \mathbb{G}_m \times \X_{\L_J} (t_2t_3 t_4), \qlb).$}

\noindent Now $\X_{\L_J} (t_2t_3 t_4)$ is a Deligne-Lusztig variety associated to a Coxeter element, and therefore the cohomology of its quotient by $U_I \cap \L_J$ is given by \cite[Corollary 2.10]{Lu}. We obtain

\centers{$  \begin{array}{r@{\ \, \simeq \, \ }l} \Hc^\bullet (\X_{W_I w_0 t_1 t_2}, \qlb)^{U_I} & \big(\Hc^\bullet(\mathbb{G}_a \times  \mathbb{G}_m \times \X_{\L_J} (t_2t_3 t_4), \qlb)\big)^{U_I \cap \L_J}\\[5pt]  & \Hc^\bullet(\mathbb{G}_a \times  (\mathbb{G}_m)^2 \times \X_{\L_I} (t_2t_3), \qlb). \end{array}$}

\item  For the piece $\X_{W_I w_0 t_4 t_3}$ we proceed as above: the minimal element $z$ of $W_K w_0 t_4 t_3$ is clearly $z = w_K w_0 t_4 t_3$ since $t_4 t_3$ is $K$-reduced. We can decompose $w$ as $w = w_1 w_2 w_3$ where $w_1 = t_1$, $w_2 = t_2 = t_2 w_2'$ and $w_3 =  t_3 t_2 t_3 t_4 $. We  observe that ${}^{z} w_1 =  t_1$, ${}^z w_2' =1$ and ${}^z w_3 = t_2 t_3$ are elements of $W_K$. In addition,  we can check by explicit computation that a piece of $\X(\mathbf{w}_1 \mathbf{w}_2 \mathbf{w}_3)$ corresponding to $(W_K z,W_K z', W_K z'')$ is empty unless $z,z'$ and $z''$ lie in the same coset. Consequently,  we can apply Theorem \ref{mainthm} to relate the cohomology of $\X_{W_K z}$ to the cohomology of $\X_{\L_K}(t_1 t_2 t_3)$ and then use \cite{Lu} to obtain

\centers{$  \begin{array}{r@{\ \, \simeq \, \ }l} \Hc^\bullet (\X_{W_I w_0 t_4 t_3}, \qlb)^{U_I} & \big(\Hc^\bullet(\mathbb{G}_a \times  \mathbb{G}_m \times \X_{\L_K} (t_1 t_2t_3), \qlb)\big)^{U_I \cap \L_K}\\[5pt]  & \Hc^\bullet(\mathbb{G}_a \times  (\mathbb{G}_m)^2 \times \X_{\L_I} (t_2t_3), \qlb). \end{array}$}

\item For the open piece $\X_{W_I w_0}$ we can directly apply Theorem \ref{mainthm} by decomposing $w$ as $w =w_1 w_2$ with $w_1= t_1 \, (t_2 t_3 t_2 t_3)$ and $w_2 =  t_4$. We only have to check that $\X_{W_I w_0}(w) = \X_{(W_Iw_0,W_Iw_0)}(\mathbf{w}_1 \mathbf{w}_2)$ which can be done using $\mathsf{GAP}$.  We deduce 

\centers{$  \Hc^\bullet (\X_{W_I w_0}, \qlb)^{U_I} \, \simeq \, \Hc^\bullet((\mathbb{G}_m)^2 \times \X_{\L_I} (t_2t_3 t_2 t_3), \qlb).$}

\end{itemize}

\noindent Note that the variety $\X_{W_I w_0 t_1 t_2} \cup \X_{W_I w_0 t_4 t_3}$ is closed in $\X(w)$. Furthermore, the elements in $W_I w_0 t_1 t_2$ and $W_I w_0 t_4 t_3$ are not comparable in the Bruhat order and therefore both $\X_{W_I w_0 t_1 t_2}$ and $\X_{W_I w_0 t_4 t_3}$ are closed subvarieties of the union. In particular

\centers{$ \Hc^\bullet\big( \X_{W_I w_0 t_1 t_2} \cup \X_{W_I w_0 t_4 t_3},\qlb\big)^{U_I} \, \simeq \,   \Big(\Hc^\bullet\big(\mathbb{G}_a \times (\mathbb{G}_m)^2\times \X_{\L_I}(t_3 t_2 ), \qlb \big) \Big)^{\oplus 2}.$}

\sk

The Weyl group of $\L_I$ has type $B_2$. Let us denote by $\varepsilon$ the sign representation of $W_I$, by $\theta$  the one-dimensional representation such that $\theta(t_2) = 1$ and $\theta(t_3)=-1$ and by $r$ the reflection representation. Then the unipotent characters of $L_I$ are  $\{\mathrm{Id}, \mathrm{St},\rho_{\theta}, \rho_{\theta \varepsilon}, \rho_r, \theta_{10}\}$ where $\theta_{10}$ is the unique unipotent cuspidal character. Using \cite[Theorem 4.3.4]{DMR} we obtain

\centers{$\begin{array}{rcl} \H_{ U_I \backslash \X_{W_I w_0 }} \hskip -2mm & = & \hskip-2mm (h^2 t  + h)^2 \big(h^4 \mathrm{St} + h^5 t^2(\rho_{\theta} + \rho_{\theta \varepsilon} + 2 \theta_{10}) + h^8 t^4 \mathrm{Id} \big) 
\\[7pt]
\hskip -2mm & = & \hskip-2mm  h^6\mathrm{St} + h^7\big(2t\mathrm{St} + t^2(\rho_{\theta} + \rho_{\theta \varepsilon} + 2 \theta_{10}) \big) + h^8\big(t^2 \mathrm{St} + 2t^3(\rho_{\theta} + \rho_{\theta \varepsilon} + 2 \theta_{10}) \big) \\[4pt] &  & \hphantom{aaaaaaaiaa} + h^9 t^4 ( \rho_{\theta} + \rho_{\theta \varepsilon} + 2 \theta_{10}) + h^{10}t^4 \mathrm{Id} + 2h^{11}t^5 \mathrm{Id} +  h^{12} t^6 \mathrm{Id} \end{array} $}

\noindent and also

\centers{$\begin{array}{rcl} \H_{U_I \backslash \X_{W_I w_0 t_1 t_2}} \hskip -2mm & = & \hskip-2mm h^2 t (h^2 t  + h)^2 \big(h^2 (\mathrm{St} + t\theta_{10}) + h^3 t \rho_r + h^4 t^2 \mathrm{Id} \big) 
\\[7pt]
\hskip -2mm & = & \hskip-2mm  h^6(t \mathrm{St} + t^2 \theta_{10}) + h^7\big( t^2 (2\mathrm{St} + \rho_r) + 2 t^3 \theta_{10} \big) \\[4pt] && \hphantom{aaiaai}  + h^8\big(t^3( \mathrm{St} + 2\rho_r+ \mathrm{Id}) + t^4 \theta_{10}\big)    + h^9 t^4( \rho_r + 2 \mathrm{Id}) + h^{10} t^5 \mathrm{Id}.\end{array} $}

We observe that the unipotent characters  $\rho_\theta$, $\rho_{\theta \varepsilon}$ and $\rho_r$ appear in the cohomology of only one of the two varieties. Using the long exact sequence associated to the decomposition $U_I \backslash \X(w) \, = \, U_I \backslash \X_{W_I w_0} \cup \big( U_I \backslash \X_{W_I w_0 t_1 t_2} \cup U_I \backslash \X_{W_I w_0 t_4 t_3} \big) $ we deduce the isotypic part of these characters in the cohomology of $U_I \backslash \X(w)$. It is given by
\begin{equation}\label{isotyp1} h^7t^2(\rho_\theta + \rho_{\theta \varepsilon} + 2 \rho_r ) + h^8 t^3(2 \rho_\theta + 2 \rho_{\theta \varepsilon} + 4 \rho_r) + h^9 t^4(\rho_\theta + \rho_{\theta \varepsilon} + 2 \rho_r ).\end{equation}
\noindent The isotypic parts for the unipotent characters  $\mathrm{St}$ and $\mathrm{Id}$ fit into the following exact sequences:

\centers{$ \begin{array}{c} 0 \longrightarrow \mathrm{St} \longrightarrow \Hc^6\big(U_I \backslash \X(w)\big)_{\mathrm{St}} \longrightarrow 2 t \mathrm{St} \longrightarrow  2 t \mathrm{St} \longrightarrow \Hc^7\big(U_I \backslash \X(w) \big)_{\mathrm{St}}  \\[6pt] \longrightarrow 4t^2 \mathrm{St} \longrightarrow t^2 \mathrm{St}  \longrightarrow \Hc^8\big(U_I \backslash \X(w) \big)_{\mathrm{St}}   \longrightarrow 2t^3 \mathrm{St} \longrightarrow 0  \end{array}$}

\centers{$ 0 \longrightarrow  \Hc^8\big(U_I \backslash \X(w) \big)_{\mathrm{Id}} \longrightarrow 2 t^3 \mathrm{Id} \longrightarrow 0$}

\centers{$ \begin{array}{c} 0 \longrightarrow  \Hc^9\big(U_I \backslash \X(w) \big)_{\mathrm{Id}} 
\longrightarrow 4t^4 \mathrm{Id} \longrightarrow t^4 \mathrm{Id} \longrightarrow  \Hc^{10}\big(U_I \backslash \X(w) \big)_{\mathrm{Id}} \\[6pt] \longrightarrow 2 t^5 \mathrm{Id} \longrightarrow 2t^5 \mathrm{Id} \longrightarrow  \Hc^{11}\big(U_I \backslash \X(w) \big)_{\mathrm{Id}} \longrightarrow 0 \end{array}$}
   
 \centers{$    0 \longrightarrow    t^6 \mathrm{Id} \longrightarrow  \Hc^{12}\big(U_I \backslash \X(w) \big)_{\mathrm{Id}}\longrightarrow 0$}
   
\noindent Any morphism above is  $F$-equivariant so that we can consider each power of $t$ separately. On the other hand, the only unipotent character of $G$ whose restriction is $\mathrm{St}_{L_I}$ (resp. $\mathrm{Id}_{L_I}$) is $\mathrm{St}_G$ (resp. $\mathrm{Id}_G$).  But from  \cite[Proposition 3.3.14 and 3.3.15]{DMR} we know exactly where these characters can appear in the cohomology of $\X(w)$ as well as the corresponding eigenvalue of $F$. Using \ref{isotyp1} we deduce that $t \mathrm{St}$ (resp. $t^2 \mathrm{St}$) cannot appear in $\Hc^6\big(U_I \backslash \X(w)\big)$ or in $\Hc^7\big(U_I \backslash \X(w)\big)$ (resp. in $\Hc^8\big(U_I \backslash \X(w)\big)$) and that  $t^4 \mathrm{Id}$ (resp. $t^5 \mathrm{Id}$) cannot appear in  $\Hc^{10}\big(U_I \backslash \X(w)\big)$ (resp.  in $\Hc^{11}\big(U_I \backslash \X(w)\big)$). With the previous exact sequences, this forces the isotypic part of $\mathrm{St}$ and $\mathrm{Id}$ in the cohomology of $U_I \backslash \X(w)$ to be

\centers{$h^6 \mathrm{St} +  3h^7t^2 \mathrm{St} + h^8 t^3 (2 \mathrm{St} + 2 \mathrm{Id})  + 3 h^9 t^4 \mathrm{Id} + h^{12} t^6 \mathrm{Id}.$}

\sk

Together with \ref{isotyp1} we finally obtain

\begin{prop}\label{5prop4}Let $w = t_1 t_2 t_3  t_2 t_3 t_4 $ and $I = \{t_2,t_3\}$. The characters of the principal series in the cohomology of $U_I \backslash \X(w)$ are given by

\centers{$\begin{array}{c} h^6 \mathrm{St} +  h^7t^2(3 \mathrm{St} + \rho_\theta + \rho_{\theta \varepsilon} + 2 \rho_r ) + h^8 t^3(2 \mathrm{St} + 2 \rho_\theta + 2 \rho_{\theta \varepsilon} + 4 \rho_r + 2 \mathrm{Id})  \\[5pt]\hphantom{aaaaaaaaaaaaaaaaaaiaiaa} + h^9 t^4(\rho_\theta + \rho_{\theta \varepsilon} + 2 \rho_r + 3 \mathrm{Id}) + h^{12} t^6 \mathrm{Id}. \end{array}$}

\end{prop}

\begin{rmk}\label{5rem3}The long exact sequence coming from the decomposition of the variety $U_I \backslash \X(w)$ does not give  enough information to determine the $\theta_{10}$-isotypic part:

\centers{$ \begin{array}{c} 0 \longrightarrow \Hc^6\big(U_I \backslash \X(w)\big)_{\theta_{10}} \longrightarrow 2 t^2 \theta_{10} \longrightarrow  2 t^2\theta_{10}  \longrightarrow \Hc^7\big(U_I \backslash \X(w) \big)_{\theta_{10} }  \longrightarrow 4t^3 \theta_{10}  \\[6pt] \longrightarrow 4t^3 \theta_{10}  \longrightarrow \Hc^8\big(U_I \backslash \X(w) \big)_{\theta_{10} }   \longrightarrow 2t^4 \theta_{10}  \longrightarrow 2t^4 \theta_{10} \longrightarrow \Hc^9\big(U_I \backslash \X(w) \big)_{\theta_{10} }   \longrightarrow 0.  \end{array}$}

\noindent One could nonetheless hope that in this particular situation the boundary maps are isomorphisms, which would imply that $\theta_{10}$ cannot appear in the cohomology of $U_I \backslash \X(w)$. This will be the case if and only if the graded endomorphism ring $\mathrm{End}_{G}(\Hc^\bullet(\X(w), \qlb))$ is concentrated in degree $0$. 
\end{rmk}

\mk

\noindent \textbf{\thesubsection.2.  Cuspidal characters.}\label{212} From \cite{BMM} we know that the irreducible constituents of the alternating sum of the comology of $\X(w)$ are the unipotent characters in the principal $\Phi_8$-block, namely $\{\mathrm{Id}_G, \mathrm{St}_G, \phi_{9,10}, \phi_{16,5}, \phi_{9,2}\}$ for the principal series and $\{\mathrm{F}_4[-1 ], \mathrm{F}_4[\mathrm{i} ], \mathrm{F}_4[-\mathrm{i} ]\}$ for the cuspidal characters (with the notation in \cite{Car}). We observe that the restriction of these characters to $L_I$ are exactly the one obtained in the previous proposition. Since the Harish-Chandra restriction preserves the Harish-Chandra series, we can deduce the contribution of the principal series to the cohomology of $\X(w)$. The missing ones are either in the series associated to $\theta_{10}$ $-$ which we could not determine $-$ or are cuspidal characters. We shall deduce the contribution of  the latter  using the results in Section \ref{5se22}. 

\sk

Recall that $G$ has $7$ cuspidal unipotent characters, namely $\mathrm{F}_4[-1 ]$, $\mathrm{F}_4[\mathrm{i} ]$, $\mathrm{F}_4[-\mathrm{i} ]$, $\mathrm{F}_4[\theta ]$, $\mathrm{F}_4[\theta^2 ]$, $\mathrm{F}_4^\mathrm{I}[1 ]$ and $\mathrm{F}_4^{\mathrm{II}}[1 ]$ where $\mathrm{i}$ (resp. $\theta$) is a primitive $4$th root of unity (resp. a primitive $3$rd root of unity). Let $\rho$ be a cuspidal unipotent character and let $v \leq w$. By cuspidality $\rho$ cannot appear in the cohomology of Deligne-Lusztig varieties associated to elements lying in a proper parabolic subgroup of $W$. In particular it cannot appear in the cohomology of $\X(v)$ or $\X(\underline{v})$ unless $v$ is in the following set

\centers{$\mathcal{V} = \big\{ w, t_1 t_2 t_3 t_2 t_4, t_1 t_3 t_2 t_3 t_4, t_1 t_2 t_3 t_4, t_1 t_3 t_2 t_4 \big\}.$}

Define $\Z =  \X(\underline{ t_1 t_2 t_3 t_2 t_4}) \cup \X(\underline{ t_1 t_3 t_2 t_3 t_4})$ and $\Z' =   \X(\underline{t_1 t_2 t_3 t_4}) \cup \X(\underline{t_1 t_3 t_2 t_4}) $. The property \hyperref[cusp2]{$\mathrm{(C2)}$} yields the following exact sequences:

\begin{equation} \cdots \longrightarrow \ \Hc^i\big(\X(w)\big)_\rho \ \longrightarrow \Hc^i\big(\X(\underline{w})\big) _\rho \ \longrightarrow \ \Hc^i\big( \Z \big)_\rho \ \longrightarrow \cdots  \label{5eq4}
\end{equation}
\begin{equation} \cdots \longrightarrow \ \Hc^i\big(\X(t_1 t_2 t_3 t_2 t_4 )\big)_\rho \ \longrightarrow \Hc^i\big(\X(\underline{t_1 t_2 t_3 t_2 t_4 }) \big)_\rho \ \longrightarrow  \ \Hc^i\big( \Z' \big)_\rho \ \longrightarrow \cdots \label{5eq5}\end{equation}
\begin{equation} \cdots \longrightarrow \ \Hc^i\big(\X(t_1 t_3 t_2 t_3 t_4 )\big)_\rho \ \longrightarrow \Hc^i\big(\X(\underline{t_1 t_3 t_2 t_3  t_4 }) \big)_\rho \ \longrightarrow  \ \Hc^i\big( \Z' \big)_\rho \ \longrightarrow \cdots \label{5eq6}\end{equation}

\noindent Moreover, one can check that each of these compactifications is actually rationally smooth, and therefore one can use \hyperref[cusp3]{$\mathrm{(C3)}$} to compute the cuspidal part of their cohomology, denoted by $\underline{\mathrm{H}}_{\X}(t^{1/2},h)$. They are given by
\begin{equation} \underline{\H}_{\X(\underline{w})} \, = \, h^6 t^3 \big( \mathrm{F}_4[-1 ] + \mathrm{F}_4[\mathrm{i} ] +\mathrm{F}_4[-\mathrm{i} ] + 2 \mathrm{F}_4[ \theta ] + 2 \mathrm{F}_4[\theta^2 ] \big) \label{5eq7} \end{equation}
\leftcenters{and}{$ \ \  \underline{\H}_{\X(\underline{t_1 t_2 t_3 t_2 t_4} )}\, = \, \underline{\H}_{\X(\underline{t_1 t_3 t_2 t_3 t_4} )} \, = \, (h^4 t^2 + h^6 t^3) \big( \mathrm{F}_4[\mathrm{i} ] +\mathrm{F}_4[-\mathrm{i} ] +  \mathrm{F}_4[ \theta ] + \mathrm{F}_4[\theta^2 ] \big).  $}

\noindent Furthermore, the elements  $t_1 t_2 t_3 t_4$ and $t_1 t_3 t_2 t_4$ are minimal in the set $\mathcal{V}$ for the Bruhat order, so that for any unipotent cuspidal character $\rho$

\centers{$ \Hc^i(\Z')_\rho \, \simeq \, \Hc^i\big(\X(t_1 t_2 t_3 t_4)\big)_\rho \oplus \Hc^i\big(\X(t_1 t_3 t_2 t_4)\big)_\rho \, \simeq \, \Hc^i\big(\X(c)\big)_\rho^{\oplus 2} $}

\noindent where $c$ is any Coxeter element of $W$. Using \cite[table 7.3]{Lu} we deduce that 

\centers{$ \underline{\H}_{\Z'} \, = \, 2h^4 t^2 \big(  \mathrm{F}_4[\mathrm{i} ] +\mathrm{F}_4[-\mathrm{i} ] +  \mathrm{F}_4[ \theta ] + \mathrm{F}_4[\theta^2 ] \big). $}

\noindent Together with \ref{5eq5} and \ref{5eq6}, and the fact that the cohomology of   $\X(t_1 t_2 t_3 t_2 t_4)$ and $\X(t_1 t_3  t_2 t_3 t_4)$ vanishes in degree $4$, we obtain

\centers{$ \underline{\H}_{\X(t_1 t_2 t_3 t_2 t_4)} \, = \, \underline{\H}_{\X(t_1 t_3 t_2 t_3 t_4)} \, = \, (h^5 t^2  + h^6 t^3)  \big( \mathrm{F}_4[\mathrm{i} ] +\mathrm{F}_4[-\mathrm{i} ] +  \mathrm{F}_4[ \theta ] + \mathrm{F}_4[\theta^2 ] \big).   $}

From these results, one can now partially determine the cohomology of $\Z$: for any unipotent cuspidal character, we use the following exact sequence

\centers{ $\cdots \longrightarrow \Hc^i\big(\X(t_1 t_2 t_3 t_2 t_4 )\big)_\rho \oplus \Hc^i\big(\X(t_1 t_3 t_2 t_3 t_4)\big)_\rho \longrightarrow \Hc^i\big(\Z \big)_\rho \longrightarrow  \Hc^i\big( \Z' \big)_\rho  \longrightarrow \cdots $}

\noindent to deduce that there exist integers $0 \leq \varepsilon_i \leq 2$ such that

\centers{$ \begin{array}{rcl} \underline{\H}_{\Z} & \hskip-2mm  =  & \hskip -2mm  (h^4 + h^5) t^2 \big(  \varepsilon_1 \mathrm{F}_4[\mathrm{i} ] +\varepsilon_2 \mathrm{F}_4[-\mathrm{i} ] +  \varepsilon_3 \mathrm{F}_4[ \theta ] + \varepsilon_4 \mathrm{F}_4[\theta^2 ]  \big) \\[5pt] & & + \,  2h^6 t^3 \big(  \mathrm{F}_4[\mathrm{i} ] +\mathrm{F}_4[-\mathrm{i} ] +  \mathrm{F}_4[ \theta ] + \mathrm{F}_4[\theta^2 ] \big). \end{array}$}

\noindent However  \ref{5eq4} forces each character $\varepsilon_i \rho$ to be a component of  $\Hc^5\big(\X(w)\big)$ since $\Hc^4\big(\X(\underline{w})\big)$ is zero by \ref{5eq7}. But the cohomology of  $\X(w)$ vanishes outside the degrees $6, \ldots, 12$, and hence the $\varepsilon_i$'s must be zero. Consequently, the exact sequence \ref{5eq4} can be decomposed into

\centers{$ 0 \longrightarrow \Hc^6\big(\X(w) \big)_{\mathrm{F}_4[-1]} \longrightarrow t^3  \mathrm{F}_4[-1] \longrightarrow 0$}

\centers{$ 0 \longrightarrow \Hc^6\big(\X(w) \big)_{\mathrm{F}_4[\pm\mathrm{i}]} \longrightarrow t^3  \mathrm{F}_4[\pm\mathrm{i}] \longrightarrow  2 t^3\mathrm{F}_4[ \pm\mathrm{i}] \longrightarrow \Hc^7\big(\X(w) \big)_{\mathrm{F}_4[\pm\mathrm{i}]}  \longrightarrow 0$ }



\centers{$ 0 \longrightarrow \Hc^6\big(\X(w) \big)_{\mathrm{F}_4[\theta^j]} \longrightarrow 2t^3  \mathrm{F}_4[\theta^j] \longrightarrow  2 t^3\mathrm{F}_4[ \theta^j] \longrightarrow \Hc^7\big(\X(w) \big)_{\mathrm{F}_4[\theta^j]}  \longrightarrow 0.$ }

\sk

We use \hyperref[cusp4]{$\mathrm{(C4)}$} to conclude: the characters $\mathrm{F}_4[\pm\mathrm{i}]$,  and $\mathrm{F}_4[\theta^j]$ already occur in the cohomology of the Deligne-Lusztig variety associated to a Coxeter element. Since $\mathbf{w}$ is not $F$-conjugate to a Coxeter element, they cannot appear in $\Hc^6 \big(\X(w)\big)$ with an eigenvalue of absolute value $q^3$, and the previous exact sequences determine the cuspidal part of the cohomology of $\X(w)$.

\begin{prop}\label{5prop5}Let $w = t_1 t_2 t_3 t_2 t_4$. The cuspidal part of the  cohomology of $\X(w)$ is given by

\centers{$ h^6 t^3 \mathrm{F}_4[-1] + h^7 t^3 \big(\mathrm{F}_4[\mathrm{i}]
 + \mathrm{F}_4[-\mathrm{i}]). $}
 
 \end{prop} 

\mk

\noindent \textbf{\thesubsection.3. Cohomology of $\X(w)$.} The unipotent characters in the principal  $\Phi_8$-block $b$ are given by $ b_{\mathrm{uni}}\, = \, \big\{\mathrm{Id}, \mathrm{St}_G, \phi_{9,10}, \phi_{16,5}, \phi_{9,2}, \mathrm{F}_4[-1 ], \mathrm{F}_4[\mathrm{i} ], \mathrm{F}_4[-\mathrm{i} ]\big\}.$ From Proposition \ref{5prop4} and \ref{5prop5} we deduce the contribution to the cohomology of $\X(w)$ of any unipotent character in the block:

\begin{thm}\label{5thm1}Let $(\G,F)$ be a split group of type F$_4$ and $w$  be a good $8$-regular element. The contribution to the cohomology of the Deligne-Lusztig $\X(w)$ of the principal series and the cuspidal characters coincides with the contribution of the principal $\Phi_8$-block, and it is given by

\centers{\begin{tabular}{@{{\vrule width 1pt}\ \, }c@{\ \,{\vrule width 1pt}\ \,}c|c|c|c|c|c|c@{\, \ {\vrule width 1pt}}} \hlinewd{1pt}
$ i $ & $\vphantom{\Big(}6$ & $7$ &$ 8$ &$ 9$ & $10$ & $11$ & $12$  \\\hlinewd{1pt} $b\H^i\big(\X(w),\qlb\big) $& $\vphantom{\Big(^y} \mathrm{St} $& $q^2 \phi_{9,10} $& 
$q^3 \phi_{16,5} $& $q^4 \phi_{9,2} $& & & $q^6\mathrm{Id} $\\[5pt]
& $-q^{3} F_4[-1] $&  & & &  & & \\[5pt]
& &$\mathrm{i} q^{3} F_4[\mathrm{i}] $& & &  & & \\[5pt]
& & $-\mathrm{i} q^{3} F_4[-\mathrm{i}] $&  & & & & \\[5pt] \hlinewd{1pt}
\end{tabular}}

\end{thm}

\subsection{\texorpdfstring{$9$}{9}-regular elements for groups of type \texorpdfstring{$E_6$}{E6}} 

In this section we assume that $(\G,F)$ is a split group of type $E_6$. The largest regular number (excluding the Coxeter number) being $9$, we are interested in computing the cohomology of $\X(w)$ for any $9$th root of $\boldsymbol \pi$, or equivalently for any good $9$-regular element. We will label the simple reflections as follows

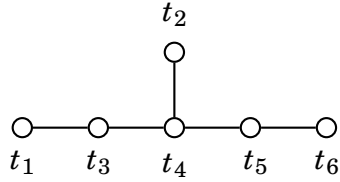
\centers{\begin{pspicture}(4,2.1)

  \cnode(0,0.5){4pt}{A}
  \cnode(2,1.5){4pt}{G}
  \cnode(1,0.5){4pt}{B}
  \cnode(2,0.5){4pt}{C}
  \cnode(3,0.5){4pt}{D}
  \cnode(4,0.5){4pt}{E}

  \ncline[nodesep=0pt]{A}{B}\nbput[npos=-0.2]{$\vphantom{\big(}t_1$}\nbput[npos=1.2]{$\vphantom{\big(}t_3$}
  \ncline[nodesep=0pt]{G}{C}\ncput[npos=-0.89]{$t_2$}\ncput[npos=1.89]{$t_4$}
  \ncline[nodesep=0pt]{B}{C}
  \ncline[nodesep=0pt]{C}{D}
  \ncline[nodesep=0pt]{D}{E}\nbput[npos=-0.1]{$\vphantom{\big(}t_5$}\nbput[npos=1.2]{$\vphantom{\big(}t_6$}
\end{pspicture}}

As before, we may, and we will, consider a particular root of $\boldsymbol \pi$, since the cohomology of the corresponding Deligne-Lusztig variety $\X(w)$ does not depend on this choice. We set

\centers{$ w =t_1 t_3  t_4 t_3 t_2 t_4  t_5  t_6. $}

\mk

\noindent \textbf{\thesubsection.1. Cohomology of $U_I \backslash \X(w)$.} We decompose the quotient of $\X(w)$ by $U_I$ for $I = \{t_2,t_3,t_4,t_5\}$. The situation is similar to the one studied in Section \hyperref[211]{\ref*{211}.1}: a piece $\X_x$ is non-empty if and only if $W_I x$ is one of the three cosets among $W_I w_0$, $W_I w_0 t_6 t_5 t_4$ and $W_I w_0 t_1 t_3$.

\begin{itemize}

\item Let $J = S \smallsetminus \{t_1\}$. We have $W_J w_0 t_1 t_3 = W_J w_0$ and therefore the piece corresponding to $W_I w_0 t_6 t_5 t_4$ is stable by the action of $P_J$.  Let $y$ be the minimal element of $W_J w_0 t_6 t_5 t_4$. Since ${}^{w_0} (t_6 t_5 t_4) = t_1 t_3 t_4$ is $J$-reduced, $y = w_J w_0 t_6 t_5 t_4$. Let us decompose $w$ as $w= w_1 w_2 w_3$ with $w_1 = t_1$, $w_2 =  t_3 = t_3 w_2'$ and $w_3 = t_4 t_3 t_2 t_4 t_5 t_6$. Then ${}^{y} w_1= t_6$, ${}^y w_2'=1$ and ${}^{y} w_3 = t_3 t_5 t_4 t_2$ are all elements of $W_J$. In addition, they satisfy the assumptions of Theorem \ref{mainthm} (see also Remark \ref{rmk1}) so that 

\centers{$  \Hc^\bullet (\X_{W_J y}, \qlb)^{U_J} \, \simeq \, \Hc^\bullet(\mathbb{G}_a \times  \mathbb{G}_m \times \X_{\L_J} (t_6 t_3 t_5 t_4 t_2), \qlb).$}

\noindent Now $\X_{\L_J} (t_6 t_3 t_5 t_4 t_2)$ is a Deligne-Lusztig variety associated to a Coxeter element, and therefore the cohomology of its quotient by $U_I \cap \L_J$ is given by \cite{Lu}. We obtain

\centers{$  \begin{array}{r@{\ \, \simeq \, \ }l} \Hc^\bullet (\X_{W_I w_0 t_6 t_5 t_4}, \qlb)^{U_I} & \big(\Hc^\bullet(\mathbb{G}_a \times  \mathbb{G}_m \times \X_{\L_J} (t_6 t_3 t_5 t_4 t_2), \qlb)\big)^{U_I \cap \L_J}\\[5pt]  & \Hc^\bullet(\mathbb{G}_a \times  (\mathbb{G}_m)^2 \times \X_{\L_I} (t_2 t_5 t_4 t_3), \qlb). \end{array}$}

\item  For the piece $\X_{W_I w_0 t_1 t_3}$, we proceed as above: let $K = S \smallsetminus \{t_6\}$ and $z$ be the minimal element of $W_K w_0 t_1 t_3$. It is clearly $z = w_K w_0 t_1 t_3$ since $t_6 t_5$ is $K$-reduced. We can decompose $w$ as $w =w_1 w_2$ where $w_1 = t_1 t_3 t_4 t_3 t_2$ and $w_2 = t_4 (t_5 t_6) = t_4 w_2'$. We have ${}^{z} w_1 = t_2 t_4 t_5 $ and ${}^z w_2' = t_3 t_1$ so that with Theorem \ref{mainthm} and \cite{Lu} we obtain

\centers{$  \begin{array}{r@{\ \, \simeq \, \ }l} \Hc^\bullet (\X_{W_I w_0 t_1 t_3}, \qlb)^{U_I} & \big(\Hc^\bullet(\mathbb{G}_a \times  \mathbb{G}_m \times \X_{\L_K} (t_2 t_4t_5 t_3 t_1), \qlb)\big)^{U_I \cap \L_K}\\[5pt]  & \Hc^\bullet(\mathbb{G}_a \times  (\mathbb{G}_m)^2 \times \X_{\L_I} (t_5 t_4 t_2 t_3), \qlb). \end{array}$}

\item For the open piece $\X_{W_I w_0}$ we can directly apply Theorem \ref{mainthm} by decomposing $w$ as $w = w_1 w_2 $ with $w_1 = t_1 \, (t_3 t_4 t_3 t_2 t_4 t_5)$ and $w_2 = t_6$. We deduce 

\centers{$  \Hc^\bullet (\X_{W_I w_0}, \qlb)^{U_I} \, \simeq \, \Hc^\bullet((\mathbb{G}_m)^2 \times \X_{\L_I} (t_5t_4 t_5 t_2 t_4 t_3), \qlb).$}

\end{itemize}

\noindent By the properties of the Bruhat order the varieties $ \X_{W_I w_0 t_6 t_5 t_4}$ and $\X_{W_I w_0 t_1 t_3}$ are both closed subvarieties of $\X(w)$. Therefore the cohomology of the union $ \X_f =  U_I \backslash \X_{W_I w_0 t_6 t_5 t_4} \cup U_I \backslash \X_{W_I w_0 t_1 t_3}$ can be deduced from \cite[Table 7.3]{Lu}  whereas the cohomology of $\X_o = U_I \backslash \X_{W_I w_0}$ is given by  \cite[Theorem 12.4]{DM2}:

\centers{$\begin{array}{rcl} \H_{\X_o} \hskip -2mm & = & \hskip-2mm (h^2 t  + h)^2 \big(h^6 \mathrm{St} + h^7\big( t^2(\rho_{1^2 +} + \rho_{1^2 -} + \rho_{21^2}) + 2t^3 \mathrm{D}_4 \big) + 2h^8 t^3 \rho_{1.21} \\[4pt] & & \phantom{ aaaaaaaaaaaaaaaaaaaaaaaaaaaa} + h^9 t^4 (\rho_{2+} + \rho_{2-} + \rho_{31} ) + h^{12} t^6 \mathrm{Id} \big) 
 \\[12pt]
 \H_{\X_f} \hskip -2mm & = & \hskip-2mm 2 h^2 t (h^2 t  + h)^2 \big(h^4 (\mathrm{St} + t^2\mathrm{D}_4) + h^5 t \rho_{1.1^3} + h^6 t^2 \rho_{1^2.2} + h^7 t^3 \rho_{1.3} + h^8 t^4\mathrm{Id} \big) 
\end{array} $}

\noindent where $\rho_\lambda$ is the unipotent character (in the principal series) associated to the character $\lambda$ of $W_I$ and $\mathrm{D}_4$ is the unique unipotent cuspidal character of $L_I$. 

\sk

As before, any character in the principal series which is different from $\mathrm{St}$ and $\mathrm{Id}$ cannot appear in the cohomology of both of the varieties, so that the isotypic part on the cohomology of $U_I \backslash \X(w)$ is the sum of the isotypic part on $\Hc^\bullet(\X_f)$ and $\Hc^\bullet(\X_o)$. For the characters  $\mathrm{St}$ and $\mathrm{Id}$, we proceed exactly as in Section \hyperref[211]{\ref*{211}.1} using  \cite[Proposition 3.3.14]{DMR} and  \cite[Proposition 3.3.15]{DMR}.

\begin{prop}\label{5prop6}Let $w = t_1 t_3  t_4 t_3 t_2 t_4  t_5  t_6 $ and $I = \{t_2,t_3,t_4,t_5\}$. The contribution of the characters in the principal series to the cohomology of $U_I \backslash \X(w)$ is given by

\centers{$ \begin{array}{l} \hphantom{+ \,}  h^8 \mathrm{St}+ h^9 t^2 \big( 3 \mathrm{St} + \rho_{1^2 +} + \rho_{1^2 -} + \rho_{21^2} + 2\rho_{1.1^3} \big) \\[5pt] + \, h^{10} t^3 \big( 2 \mathrm{St} +  2\rho_{1^2 +} + 2\rho_{1^2 -} + 2\rho_{21^2} + 4\rho_{1.1^3} + 2 \rho_{1.21} + 2 \rho_{1^2.2} \big) \\[5pt]
+ \, h^{11} t^4 \big(\rho_{1^2 +} + \rho_{1^2 -} + \rho_{21^2} + 2\rho_{1.1^3} + 4 \rho_{1.21} + 4 \rho_{1^2.2} + \rho_{2+} + \rho_{2-} + \rho_{31} + 2\rho_{1.3}\big)
\\[5pt]  + \, h^{12} t^5 \big( 2 \rho_{1.21} + 2 \rho_{1^2.2} + 2\rho_{2+} + 2\rho_{2-} + 2\rho_{31} + 4\rho_{1.3} + 2 \mathrm{Id} \big) \\[5pt]
+ \, h^{13} t^6 \big(\rho_{2+} + \rho_{2-} + \rho_{31} + 2\rho_{1.3} + 3\mathrm{Id}\big) + h^{16} t^8 \mathrm{Id}.
\end{array}$}

\end{prop}

\begin{rmk}\label{5rem4}Unfortunately, this method is not sufficient for determining the $\mathrm{D}_4$-isotypic part (see also Remark \ref{5rem3}).
\end{rmk}

\noindent \textbf{\thesubsection.2. Cuspidal characters.} The group $G$ has only two cuspidal characters, denoted by $\mathrm{E}_6[\theta]$ and $\mathrm{E}_6[\theta^2]$ where $\theta$ is a primitive $3$rd root of unity. In order to determine they contribution to the cohomology  of $\X(w)$ , we use the compactifications $\X(\underline{v})$ for $v \leq w$. However, unlike the type $F_4$, they are not always rationally smooth and we shall rather work with "bigger" compactifications, obtained by underlining all the simple reflections. For details on the explicit computations we refer to Section  \hyperref[212]{\ref*{212}.2}. We start by defining the following closed subvariety of $\overline{\X}(w)$: 

\centers{ $\Z = \X(\underline{t}_1 \underline{t}_4 \underline{t}_3\underline{t}_2\underline{t}_4\underline{t}_5\underline{t}_6) \cup \X(\underline{t}_1 \underline{t}_3 \underline{t}_3\underline{t}_2\underline{t}_4\underline{t}_5\underline{t}_6) \cup \X(\underline{t}_1 \underline{t}_3 \underline{t}_4 \underline{t}_2\underline{t}_4\underline{t}_5\underline{t}_6) \cup \X(\underline{t}_1 \underline{t}_3 \underline{t}_4 \underline{t}_3\underline{t}_2\underline{t}_5\underline{t}_6)$}

\noindent so that we obtain, for any cuspidal character $\rho$, a long exact sequence  
\begin{equation} \cdots \longrightarrow \ \Hc^i\big(\X(w)\big)_\rho \ \longrightarrow \Hc^i\big(\overline{\X}(w)\big) _\rho \ \longrightarrow \ \Hc^i\big(\Z \big)_\rho \ \longrightarrow \cdots  \label{5eq8}\end{equation}
\noindent We determine the cuspidal part of $\Z$ as follows: we compute, for any element $v \in \{t_1 t_4 t_3 t_2 t_4 t_5 t_6, t_1 t_3 t_4 t_2 t_4 t_5 t_6,t_1 t_3 t_4 t_3 t_2  t_5 t_6\}$

\centers{${ \underline{\H}}_{\X(\mathbf{v})} \, = \, { \underline{\H}}_{\X(v)} \, = \, (h^7 t^3 + h^8 t^4) \big( \mathrm{E}_6[\theta] +  \mathrm{E}_6[\theta^2] \big)$}

\noindent by means of the following exact sequences

\centers{$ \cdots \longrightarrow \ \Hc^i\big(\X(v)\big)_\rho \ \longrightarrow \Hc^i\big(\overline{\X}(v)\big) _\rho \ \longrightarrow \ \big(\Hc^i( \overline{\X}(c))_\rho\big)^{\oplus 2}  \ \longrightarrow \cdots$}

\noindent and the precise values

\leftcenters{and}{$\begin{array}[b]{r@{\, \ = \, \ }l} \underline{\H}_{\overline{\X}(v)} & (h^6 t^3 + h^8 t^4) \big( \mathrm{E}_6[\theta] +  \mathrm{E}_6[\theta^2] \big) \\[5pt] \underline{\H}_{\overline{\X}(c)} & h^6 t^3 \big( \mathrm{E}_6[\theta] +  \mathrm{E}_6[\theta^2] \big) \end{array}$} 

\noindent that can be found using  \hyperref[cusp3]{$\mathrm{(C3)}$}. Note that we have also used the fact that the cohomology of  $\X(v)$ is zero outside the degrees $7,\ldots,14$. For the element  $\mathbf{v} = \mathbf{t}_1\mathbf{t}_3\mathbf{t}_3\mathbf{t}_2\mathbf{t}_4\mathbf{t}_5\mathbf{t}_6$ we use \cite[Proposition 3.2.10]{DMR} and we obtain the same value again:

\centers{${ \underline{\H}}_{\X(\mathbf{v})} \, = \, (h^2t+h) \,\underline{\H}_{\overline{\X}(c)} \, = \, (h^7 t^3 + h^8 t^4) \big( \mathrm{E}_6[\theta] +  \mathrm{E}_6[\theta^2] \big)$}

\noindent In particular, the cohomology of $\Z$ fits into the following long exact sequence

\centers{$ \cdots \longrightarrow \ \big(\Hc^i(\X(\mathbf{v}))_\rho\big)^{\oplus 4}  \ \longrightarrow \Hc^i(\Z ) _\rho \ \longrightarrow \ \big(\Hc^i( \overline{\X}(c))_\rho\big)^{\oplus 4}  \ \longrightarrow \cdots$}

\noindent We claim that
\begin{equation} \underline{\H}_{\Z} \, = \, 4 h^8 t^4\big( \mathrm{E}_6[\theta] +  \mathrm{E}_6[\theta^2]\big). \label{5eq9}\end{equation}

\noindent Again, the exact sequence itself is not enough to compute this value, but it can be deduced from the following properties:

\begin{itemize}
\item the cohomology of $\X(w)$ vanishes in degree $7$ by  \hyperref[cusp1]{$\mathrm{(C1)}$};
\item $\underline{\H}_{\overline{\X}(w)} = 3 h^8 t^4 \big( \mathrm{E}_6[\theta] +  \mathrm{E}_6[\theta^2]\big)$ which forces in particular  $\Hc^6(\overline{\X}(w))$ to have no cuspidal constituent.

\end{itemize}

\noindent These properties, together with \ref{5eq8}, ensure that the coefficient of $h^6$ in $ \underline{\H}_{\Z} $ is zero, and we deduce \ref{5eq9}.

\sk

Consequently, the decomposition $\overline{\X}(w) = \X(w) \cup \Z$ yields the following exact sequence for any cuspidal character $\rho$:
\centers{$ 0 \longrightarrow \Hc^8(\X(w))_\rho \longrightarrow 3t^4 \rho \longrightarrow 4 t^4 \rho \longrightarrow \Hc^9(\X(w))_\rho \longrightarrow 0. $}

\noindent Finally, by \hyperref[cusp4]{$\mathrm{(C4)}$} the group $\Hc^8(\X(w))$ cannot contain any unipotent cuspidal character with an eigenvalue of absolute value $q^4$ and we obtain:

\begin{prop}\label{5prop7}Let $w =  t_1 t_3  t_4 t_3 t_2 t_4  t_5  t_6$. The contribution of the cuspidal characters of $G$ to the cohomology of $\X(w)$ is given by

\centers{$ h^9 t^4 \big( \mathrm{E}_6[\theta] +  \mathrm{E}_6[\theta^2]\big).$}

\end{prop}

\bk

\noindent \textbf{\thesubsection.3. Cohomology of $\X(w)$.} By \cite{BMi1}, the irreductible constituents of the virtual character associated to the cohomology of $\X(w)$ are exactly the unipotent characters in the principal $\Phi_9$-block, namely $b_{\mathrm{uni}} \, = \, \big\{ \mathrm{Id}_G, \mathrm{St}_G, \phi_{20,20},\phi_{64,13} , \phi_{90,8}, \phi_{64,4},$ $\phi_{20,2}, \mathrm{E}_6[\theta], \mathrm{E}_6[\theta^2]\big\}.$ By looking at the Harish-Chandra restriction of these characters, we can deduce from Proposition \ref{5prop6} and \ref{5prop7} the following theorem:

\begin{thm}\label{5thm2}Let $(\G,F)$ be a split group of type $E_6$ and $w$ be a good $9$-regular element of $W$. The contribution to the cohomology of the Deligne-Lusztig $\X(w)$ of the principal series and the cuspidal characters coincides with the contribution of the principal $\Phi_9$-block, and it is given by

\centers{\begin{tabular}{@{{\vrule width 1pt}\  }c@{\ {\vrule width 1pt}\ \, }c|@{\ }c@{\ }|@{\ }c@{\ }|c|c|c|@{\,}c@{\,}|@{\,}c@{\,}|c@{\ {\vrule width 1pt}}} 
\hlinewd{1pt}  
$i $&$ \vphantom{\Big(}8$ & $9$ & $10$ & $11$ &$ 12$ &$ 13$ & $14 $& $15$ & $16$ \\\hlinewd{1pt} $b\H^i(\X(w),\qlb) $&$ \vphantom{\Big(^y} \mathrm{St}$ & $q^2 \phi_{20,20} $& 
$q^3 \phi_{64,13} $&$ q^4 \phi_{90,8} $& $q^5 \phi_{64,4} $& $q^6 \phi_{20,2}$ & & &$ q^8\mathrm{Id} $\\[5pt]
& & $\theta q^{4} E_6[\theta] $& & & & & & & \\[5pt]
& & $\theta^2 q^{4} E_6[\theta^2] $& & & & & & & \\[5pt] \hlinewd{1pt}\end{tabular}}

\end{thm}

Conjecturally, for good regular elements, there should be no cancelation in the virtual character $\sum (-1)^i\Hc^i (\X(w),\qlb) \in \, K_0(G$-$\mathrm{mod})$ \cite[Conjecture 5.7]{BMi2}. In particular, the series associated to the cuspidal character of $D_4$ should not appear in the cohomology of $\X(w)$:

\begin{hyp}\label{5hyp1}For good $9$-regular elements in $E_6$, the cohomology of $\X(w)$  has no constituent in the Harish-Chandra series associated to the cuspidal representation of $D_4$.
\end{hyp}

This assumption will be essential to study the contribution of the $D_4$-series for groups of type $E_7$ and $E_8$ (see Theorem \ref{5thm3} and \ref{5thm4}).

\subsection{\texorpdfstring{$14$}{14}-regular elements for groups of type \texorpdfstring{$E_7$}{E7}}

We now assume that $(\G,F)$ is a split group of type $E_7$ and we are interested in computing the cohomology of Deligne-Lusztig varieties associated to good $14$-regular elements. We will label the simple reflections according to the following Dynkin diagram

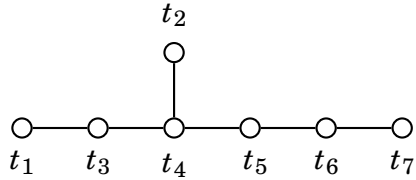
\centers{\begin{pspicture}(5,2.1)

  \cnode(0,0.5){4pt}{A}
  \cnode(2,1.5){4pt}{G}
  \cnode(1,0.5){4pt}{B}
  \cnode(2,0.5){4pt}{C}
  \cnode(3,0.5){4pt}{D}
  \cnode(4,0.5){4pt}{E}
  \cnode(5,0.5){4pt}{F}

  \ncline[nodesep=0pt]{A}{B}\nbput[npos=-0.2]{$\vphantom{\big(}t_1$}\nbput[npos=1.2]{$\vphantom{\big(}t_3$}
  \ncline[nodesep=0pt]{G}{C}\ncput[npos=-0.89]{$t_2$}\ncput[npos=1.89]{$t_4$}
  \ncline[nodesep=0pt]{B}{C}
  \ncline[nodesep=0pt]{C}{D}
  \ncline[nodesep=0pt]{D}{E}\nbput[npos=-0.1]{$\vphantom{\big(}t_5$}\nbput[npos=1.2]{$\vphantom{\big(}t_6$}
    \ncline[nodesep=0pt]{E}{F}\nbput[npos=1.2]{$\vphantom{\big(}t_7$}
\end{pspicture}}

\noindent and consider a specific $9$th root of $\boldsymbol \pi$:

\centers{$ w=t_7 t_6 t_5 t_4 t_5 t_2 t_4  t_3 t_1.$}

\mk

\noindent \textbf{\thesubsection.1. Cohomology of $U_I \backslash \X(w)$.} Let $I = S\smallsetminus \{t_7\}$. The group $\L_I$ has type $E_6$ and we can use the results in the previous section to compute the cohomology of the quotient of $\X(w)$ by $U_I$. In the decomposition of $\X(w)$ by $\P_I$-cosets in $\G/\B$, only two pieces are non-empty, with associated cosets $W_I w_0$ and $W_I w_0 t_7 t_6 t_5$. We can apply Theorem \ref{mainthm} in these two cases:

\begin{itemize} 

\item when $y = w_I w_0 t_7 t_6 t_5$ we decompose $w$ as $w=w_1w_2$ with $w_1 = t_7 t_6 t_5 t_4 t_5 t_2$ and $w_2 = t_4 \, (t_3 t_1) = t_4 w_2'$. We have ${}^y w_1 = t_1 t_3 t_4  t_2$ and ${}^y w_2' = t_5 t_6$ so that

\centers{$  \Hc^\bullet (\X_{W_I y}, \qlb)^{U_I} \, \simeq \, \Hc^\bullet(\mathbb{G}_a \times  \mathbb{G}_m \times \X_{\L_I} (t_1 t_3 t_4  t_2 t_5 t_6), \qlb).$}

\item for $x = w_I w_0$ we observe that $w = t_7 \, ( t_6 t_5 t_4 t_5 t_2 t_4  t_3 t_1) = t_7 w'$ with ${}^x w' \in W_J$ and deduce that
 
\centers{$  \Hc^\bullet (\X_{W_I w_0}, \qlb)^{U_I} \, \simeq \, \Hc^\bullet(\mathbb{G}_m \times \X_{\L_I} (t_1 t_3 t_4 t_3 t_2t_4 t_5 t_6), \qlb).$}

\end{itemize}

\noindent The cohomology of these varieties is known by Theorem \ref{5thm2} and \cite[Table 7.3]{Lu}. Recall that for any Coxeter element  $c_I$ of $W_I$, the cohomology of the corresponding variety is given by

\centers{$\hskip-2mm \begin{array}{rcl} \H_{\X_{\L_I}(c_I)} & \hskip -2mm = \hskip -2mm & h^6 \big( \mathrm{St} + t^2 \mathrm{D}_{4,\varepsilon} + t^3 \mathrm{E}_6[\theta] + t^3 \mathrm{E}_6[\theta^2] \big) + h^7 \big( t \phi_{6,25} + t^3 \mathrm{D}_{4,r} \big) \\[5pt] & &   + \, h^8 \big( t^2 \phi_{15,17} + t^4 \mathrm{D}_{4,\mathrm{Id}}\big) + h^9 t^3\phi_{20,10} + h^{10} t^4 \phi_{15,5} + h^{11} t^5 \phi_{6,1} + h^{12} t^6 \mathrm{Id}.
\end{array}$}

\noindent If we exclude $\mathrm{St}$ and $\mathrm{Id}$, none of the characters in the principal series that appear here can appear in the cohomology of $U_I \backslash \X_{W_I w_0}$. From that observation one can readily deduce the contribution of the principal series to the cohomology of $U_I \backslash \X(w)$. Note that in the case of $\mathrm{St}$ and $\mathrm{Id}$ we can proceed as in  Section \hyperref[211]{\ref*{211}.1}.

\begin{prop}\label{5prop8}Let $w = t_7 t_6 t_5 t_4 t_5 t_2 t_4  t_3 t_1 $ and $I = \Delta \smallsetminus \{t_7 \}$.  The contribution of the principal series to the cohomology of  $U_I \backslash \X(w)$ is given by

\centers{$ \begin{array}{l} \hphantom{+ \,}  h^9 \mathrm{St} + h^{10} t^2 \big(  \mathrm{St} + \phi_{6,25} + \phi_{20,20} \big) + h^{11} t^3 \big( \phi_{6,25} + \phi_{20,20} + \phi_{15,17} +  \phi_{64,13}\big) 
  \\[5pt]  + \, h^{12} t^4 \big(\phi_{15,17} +  \phi_{64,13} + \phi_{20,10} + \phi_{90,8}\big) 
+ h^{13} t^5 \big( \phi_{20,10} + \phi_{90,8} + \phi_{15,5} + \phi_{64,4} \big) 
   \\[5pt] + \, h^{14} t^6 \big( \phi_{15,5} + \phi_{64,4} + \phi_{6,1} + \phi_{20,2} \big) 
+ h^{15} t^7 \big(  \phi_{6,1} + \phi_{20,2}  + \mathrm{Id} \big)  
 + h^{18} t^9 \mathrm{Id}.
\end{array}$}

\end{prop}

The case of the Harish-Chandra series associated to the cuspidal character of $D_4$ remains undetermined unless we know the contribution of this series to the cohomology of the open part. However, in our situtation, none of these characters should appear, and the isotypic part on the cohomology of the union $U_I \backslash \X(w)$ should come from the Coxeter variety only.

\begin{prop}\label{5prop9}Assume that \ref{5hyp1} holds, and let  $w = t_7 t_6 t_5 t_4 t_5 t_2 t_4  t_3 t_1 $ and $I = \Delta \smallsetminus \{t_7 \}$. Then the contribution of the Harish-Chandra series associated to the cuspidal character of $D_4$ to the cohomology of $U_I \backslash \X(w)$ is given by

\centers{$ \begin{array}{l} \hphantom{+ \,}  h^9 t^3 \mathrm{D}_{4,\varepsilon} + h^{10} t^4 \big(\mathrm{D}_{4,\varepsilon} + \mathrm{D}_{4,r} \big) + h^{11} t^5  \big(\mathrm{D}_{4,r} + \mathrm{D}_{4,\mathrm{Id}}\big) + h^{12} t^6 \mathrm{D}_{4,\mathrm{Id}} .
\end{array}$}

\end{prop}

Finally, for the cuspidal characters $\mathrm{E}_6[\theta]$ and $\mathrm{E}_6[\theta^2]$, we have a long exact sequence

\centers{$\begin{array}{c}  0 \longrightarrow  \Hc^9\big(U_I \backslash \X(w)\big)_{\mathrm{E}_6[\theta]}  \longrightarrow t^4 \mathrm{E}_6[\theta] \longrightarrow t^4 \mathrm{E}_6[\theta] \longrightarrow  \Hc^{10}\big(U_I \backslash \X(w)\big)_{\mathrm{E}_6[\theta]} \\[5pt] \longrightarrow t^5 \mathrm{E}_6[\theta] \longrightarrow  t^5 \mathrm{E}_6[\theta] \longrightarrow \Hc^{11}\big(U_I \backslash \X(w) \big)_{\mathrm{E}_6[\theta]}  \longrightarrow 0. \end{array}$}

\noindent This is not enough to determine their contribution and we can only hope that they actually do not appear in the cohomology of $U_I \backslash \X(w)$.\mk

\noindent \textbf{\thesubsection.2. Cuspidal characters.} The group  $G$ has only two cuspidal unipotent characters, namely $\mathrm{E}_7[i]$ and $\mathrm{E}_7[-i]$, where $\mathrm{i}$ is a primitive $4$th root of unity. The method to determine their contribution to the cohomology is strictly identical to the case of $E_6$ and yields

\begin{prop}\label{5prop10}Let $w =  t_7 t_6 t_5 t_4 t_5 t_2 t_4  t_3 t_1 $. The cuspidal part of the cohomology of $\X(w)$ is given by

\centers{$ h^{10} t^{9/2} \big( \mathrm{E}_7[\mathrm{i}] +  \mathrm{E}_7[-\mathrm{i}]\big).$}

\end{prop}

\mk

\noindent \textbf{\thesubsection.3. Cohomology of $\X(w)$.} By combining Proposition \ref{5prop8} and \ref{5prop9}, we obtain the Harish-Chandra restriction to $E_6$ of the cohomology of the variety $\X(w)$. If we compare these to the restriction of the characters in the principal $\Phi_{14}$-block $b_{\mathrm{uni}}  =  \{ \mathrm{St}_G, \mathrm{Id}_G, \phi_{27,37}, \phi_{105,26}, \phi_{189,17}, \phi_{189,10}, \phi_{105,5}, \phi_{27,2}, \mathrm{D}_{4,1^3.}, \mathrm{D}_{4,1^2.1}, $ $\mathrm{D}_{4,1.2}, \mathrm{D}_{4,.3}, \mathrm{E}_7[\mathrm{i}],\mathrm{E}_7[-\mathrm{i}] \}$  (and the fact that these actually occur as constituents of the cohomology) we deduce their exact contribution. Adding the cuspidal characters obtained in  \ref{5prop10}, we get

\begin{thm}\label{5thm3}Let $(\G,F)$ be a split group of type $E_7$ and $w$ be a good $14$-regular element of $W$. The contribution to the cohomology of the Deligne-Lusztig $\X(w)$ of the principal series, the $D_4$-series and the cuspidal characters coincides with the contribution of the principal $\Phi_{14}$-block, and it is given by

\centers{ \begin{tabular}{l} \begin{tabular}{@{{\vrule width 1pt}\,\ }c@{\,\ {\vrule width 1pt}\, \ }c|c|c|c|c@{\, \ {\vrule width 1pt}}} \hlinewd{1pt} 
$i$ &$ \vphantom{\Big(} 9 $&$ 10 $& $11 $& $12 $& $13 $ 
 \\\hlinewd{1pt}$ b\H^i(\X(w),\qlb) $& $\vphantom{\Big(^y} \mathrm{St} $& $q^2 \phi_{27,37} $& $q^3 \phi_{105,26} $&$ q^4 \phi_{189,17} $&$ q^5 \phi_{189,10}$   
 \\[5pt]
&$\color{purple}-q^3 D_{4,1^3.} $& $\color{purple}-q^4 D_{4,1^2.1}$ & $\color{purple}-q^5 D_{4,1.2}$ 
& $\color{purple}-q^6 D_{4,.3} $&   \\[5pt]
& & $\mathrm{i} q^{9/2} E_7[\mathrm{i}]$ & & &  \\[5pt]
& & $-\mathrm{i} q^{9/2} E_7[-\mathrm{i}] $& & &  \\[5pt] \hlinewd{1pt}
\end{tabular}\\  \\
\begin{tabular}{@{{\vrule width 1pt}\,\ }c@{\,\ {\vrule width 1pt}\,\ }c|c|c|c|c@{\,\ {\vrule width 1pt}}} \hlinewd{1pt} $i$ & $\vphantom{\Big(} 14 $& $15 $& $16$ & $17 $& $18$ \\\hlinewd{1pt} $b\H^i(\X(w),\qlb) $& $\vphantom{\Big(^y}  q^6 \phi_{105,5} $& $q^7 \phi_{27,2} $& & & $q^9\mathrm{Id} $ \\[5pt]
\hlinewd{1pt}
\end{tabular} \end{tabular}}

\sk
\noindent where the $D_4$-series is given under Assumption \ref{5hyp1}.
\end{thm}

In our situation, the non-cancellation for the corresponding Deligne-Lusztig virtual character is equivalent to the following:

\begin{hyp}\label{5hyp2}The characters lying in the Harish-Chandra series associated to the cuspidal characters E$_6[\theta]$ and E$_6[\theta^2]$ do not appear in the cohomology of $\X(w)$. \end{hyp}

\subsection{\texorpdfstring{$24$}{24}-regular elements for groups of type \texorpdfstring{$E_8$}{E8}}

We close this section by studying the cohomology of Deligne-Lusztig varieties associated to good $24$-regular elements in $E_8$. We will label the simple reflections as follows

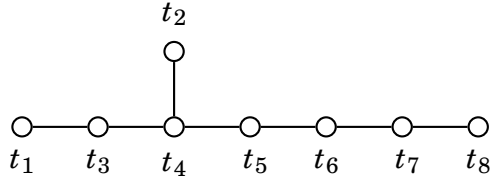
\centers{\begin{pspicture}(6,2.1)

  \cnode(0,0.5){4pt}{A}
  \cnode(2,1.5){4pt}{G}
  \cnode(1,0.5){4pt}{B}
  \cnode(2,0.5){4pt}{C}
  \cnode(3,0.5){4pt}{D}
  \cnode(4,0.5){4pt}{E}
  \cnode(5,0.5){4pt}{F}
  \cnode(6,0.5){4pt}{H}

  \ncline[nodesep=0pt]{A}{B}\nbput[npos=-0.2]{$\vphantom{\big(}t_1$}\nbput[npos=1.2]{$\vphantom{\big(}t_3$}
  \ncline[nodesep=0pt]{G}{C}\ncput[npos=-0.89]{$t_2$}\ncput[npos=1.89]{$t_4$}
  \ncline[nodesep=0pt]{B}{C}
  \ncline[nodesep=0pt]{C}{D}
  \ncline[nodesep=0pt]{D}{E}\nbput[npos=-0.1]{$\vphantom{\big(}t_5$}\nbput[npos=1.2]{$\vphantom{\big(}t_6$}
    \ncline[nodesep=0pt]{E}{F}
  \ncline[nodesep=0pt]{F}{H}\nbput[npos=-0.1]{$\vphantom{\big(}t_7$}\nbput[npos=1.2]{$\vphantom{\big(}t_8$}

\end{pspicture}}

\noindent and  choose the following $24$th root of $\boldsymbol \pi$:

\centers{$ w = t_8 t_7 t_6 t_5 t_4 t_5 t_2 t_4 t_3 t_1. $}

\mk

\noindent \textbf{\thesubsection.1. Cohomology of $U_I \backslash \X(w)$.} The situation is very similar to the case of $E_7$ so we will omit the details. When $I = S \smallsetminus \{t_8\}$, the pieces corresponding to $W_I w_0$ and $W_I w_0 t_8 t_7 t_6 t_5$ are the only non-empty pieces, and the cohomology of their quotient by  $U_I$ is given by

\centers{$ \Hc^\bullet(U_I \backslash \X_{W_I w_0 t_8 t_7 t_6 t_5},\qlb)  \, \simeq \,  \Hc^\bullet\big(\mathbb{G}_a \times \mathbb{G}_m \times \X_{\L_I}( t_7 t_6 t_5 t_4 t_2 t_3 t_1),\qlb \big) $}

\leftcenters{and}{$\Hc^\bullet(U_I \backslash \X_{W_I w_0},\qlb) \, \simeq \, \Hc^\bullet\big(\mathbb{G}_m \times \X_{\L_I}(  t_7  t_6 t_5 t_4 t_5 t_2 t_4 t_3 t_1),\qlb\big).$}

\noindent The latter has been computed in the previous section, whereas the cohomology of a Deligne-Lusztig variety associated to any Coxeter element $c_I$ of $W_I$ can be deduced from \cite[Table 7.3]{Lu}:

\centers{$\begin{array}{rcl} \H_{\X_{\L_I}(c_I)} &\hskip -2mm = \hskip -2mm & h^{7} \big(\mathrm{St} + t^2 \mathrm{D}_{4,\varepsilon} + t^3 ( \mathrm{E}_6[\theta]_{\varepsilon} + \mathrm{E}_6[\theta^2]_{\varepsilon} ) + t^{7/2} ( \mathrm{E}_7[\mathrm{i}]  +  \mathrm{E}_7[-\mathrm{i}] )\big) \\[5pt]  & & + \, h^8 \big( t\phi_{7,46} + t^3 \mathrm{D}_{4,1.1^2} + t^4 ( \mathrm{E}_6[\theta]_{\mathrm{Id}} + \mathrm{E}_6[\theta^2]_{\mathrm{Id}} ) \big) 
\\[5pt] & & + \, h^{9} \big( t^2 \phi_{21,33} + t^4 \mathrm{D}_{4, 2.1} \big) + h^{10} \big(t^3 \phi_{35,22} + t^5 \mathrm{D}_{4,\mathrm{Id}} \big) 
\\[5pt] & &  + \, h^{11} t^4 \phi_{35,13} + h^{12} t^5 \phi_{21,6} + h^{13} t^6 \phi_{7,1} + h^{14} t^7 \mathrm{Id}. \\
\end{array}$}

\noindent Together with Theorem \ref{5thm3}, this is  enough to determine the contribution of the principal series:

\begin{prop}\label{5prop11}Let $w = t_8 t_7 t_6 t_5 t_4 t_5 t_2 t_4 t_3 t_1$ and $I = \Delta \smallsetminus \{t_8 \}$. The contribution of the principal series to the cohomology of $U_I \backslash \X(w)$ is given by

\centers{$ \begin{array}{l} \hphantom{+ \,}  h^{10} \mathrm{St} + h^{11} t^2 \big(  \mathrm{St} + \phi_{7,46} + \phi_{27,37} \big) + h^{12} t^3 \big( \phi_{7,46}+ \phi_{27,37} + \phi_{21,33} +  \phi_{105,26}\big) 
  \\[5pt] + \, h^{13} t^4 \big(\phi_{21,33} +  \phi_{105,26} + \phi_{35,22} + \phi_{189,17}\big) 
+ h^{14} t^5 \big( \phi_{35,22} + \phi_{189,17} + \phi_{35,13} + \phi_{189,10} \big) 
   \\[5pt] + \, h^{15} t^6 \big( \phi_{35,13}+ \phi_{189,10} + \phi_{21,6} + \phi_{105,5} \big) 
+ h^{16} t^7 \big(  \phi_{21,6} + \phi_{105,5} + \phi_{7,1} + \phi_{27,2}\big) \\[5pt]
+ \,h^{17} t^8 \big(   \phi_{7,1} + \phi_{27,2}  + \mathrm{Id}  \big)
 + h^{20} t^{10} \mathrm{Id}.
\end{array}$}

\end{prop}

The results for the intermediate series depend whether the assumptions \ref{5hyp1} and \ref{5hyp2} are satisfied. If they hold, we can easily obtain:

\begin{prop}\label{5prop12}Let $w = t_8 t_7 t_6 t_5 t_4 t_5 t_2 t_4 t_3 t_1$ and $I = \Delta \smallsetminus \{t_8 \}$. 

\begin{enumerate} 

\item[$\mathrm{(i)}$] Under Assumption \ref{5hyp1}, the contribution of the $D_4$-series to the cohomology of $U_I \backslash \X(w)$ is given by

\centers{$ \begin{array}{l} \hphantom{+ \,}  h^{10} t^3 \big( \mathrm{D}_{4,\varepsilon} + \mathrm{D}_{4, 1^3.} \big) + h^{11} t^4 \big( \mathrm{D}_{4,\varepsilon} + \mathrm{D}_{4,1^3.} + \mathrm{D}_{4,1.1^2} + \mathrm{D}_{4,1^2.1}  \big)  
\\[5pt] + \, h^{12} t^5 \big( \mathrm{D}_{4,1.1^2} + \mathrm{D}_{4,1^2.1}  + \mathrm{D}_{4,2.1} + \mathrm{D}_{4,1.2}  \big) 
\\[5pt] + \, h^{13} t^6 \big( \mathrm{D}_{4,2.1} + \mathrm{D}_{4,1.2}  + \mathrm{D}_{4,\mathrm{Id}} + \mathrm{D}_{4,.3} \big) + h^{14} t^7 \big(  \mathrm{D}_{4,\mathrm{Id}} + \mathrm{D}_{4,.3} \big). 
\end{array}$}

\item[$\mathrm{(ii)}$]   Under Assumption \ref{5hyp2}, the contribution of the $E_6$-series to the cohomology of $U_I \backslash \X(w)$ is given by

\centers{$ h^{10} t^4 \mathrm{E}_6[\theta]_\varepsilon + h^{11} t^5\big(\mathrm{E}_6[\theta]_\varepsilon + \mathrm{E}_6[\theta]_{\mathrm{Id}} \big) + h^{12} t^6 \mathrm{E}_6[\theta]_{\mathrm{Id}}$}

\leftcenters{and}{$h^{10} t^4 \mathrm{E}_6[\theta^2]_\varepsilon + h^{11} t^5\big(\mathrm{E}_6[\theta^2]_\varepsilon + \mathrm{E}_6[\theta^2]_{\mathrm{Id}} \big) + h^{12} t^6 \mathrm{E}_6[\theta^2]_{\mathrm{Id}} .$} 

\end{enumerate}
\end{prop}

\bk

\noindent \textbf{\thesubsection.2. Cuspidal characters.} The group $G$ has several cuspidal unipotent characters, denoted in \cite{Car} by  $\mathrm{E}_8[\pm \mathrm{i}],  \mathrm{E}_8[\pm \theta], \mathrm{E}_8[\pm \theta^2],  \mathrm{E}_8^\mathrm{I}[1], \mathrm{E}_8^\mathrm{II}[1]$ and   $\mathrm{E}_8[\zeta^j]$ where $\zeta$ is a primitive $5$th root of unity  and $j = 1,2,3,4$. We proceed as in the previous cases to determine they contribution to the cohomology of $\X(w)$. However, due to the large number of cuspidal characters, the calculations are a bit more tedious. 

\sk

We start by considering the closed subvariety $\Z$ of $\overline{\X}(w)$ consiting of the union of the varieties $\X(v)$ where $v$ runs over the set

\centers{$  \big\{\underline{t}_8 \underline{t}_7 \underline{t}_6  \underline{t}_4 \underline{t}_5 \underline{t}_2 \underline{t}_4 \underline{t}_3 \underline{t}_1, \ \underline{t}_8 \underline{t}_7 \underline{t}_6 \underline{t}_5 \underline{t}_4  \underline{t}_2 \underline{t}_4 \underline{t}_3 \underline{t}_1, \ \underline{t}_8 \underline{t}_7 \underline{t}_6 \underline{t}_5 \underline{t}_4 \underline{t}_5 \underline{t}_2 \underline{t}_3 \underline{t}_1, \ \underline{t}_8 \underline{t}_7 \underline{t}_6 \underline{t}_5  \underline{t}_5 \underline{t}_2 \underline{t}_4 \underline{t}_3 \underline{t}_1\big\}.$}

\noindent The cohomology of this variety fits in the following long exact sequence, for any cuspidal character $\rho$

\begin{equation} \cdots \longrightarrow \ \Hc^i\big(\X(w)\big)_\rho \ \longrightarrow \Hc^i\big(\overline{\X}(w)\big) _\rho \ \longrightarrow \ \Hc^i\big( \Z \big)_\rho \ \longrightarrow \cdots \label{5eq10}\end{equation}

\noindent The elements of the Braid monoid obtained by un-underlying the elements $v$ will be denoted by  $\mathbf{v}_1, \mathbf{v}_2, \mathbf{v}_3$ et $\mathbf{v}_4$. Note that only $\mathbf{v}_4$ is not the canonical lift of an element of $W$. For $j=1,2,3$, the cuspidal part of the cohomology of  $\X(\mathbf{v}_j) \simeq \X(v_j)$ can be deduced from the following exact sequence

\centers{ $\cdots \longrightarrow \ \Hc^i\big(\X(v_j)\big)_\rho \ \longrightarrow \Hc^i\big(\overline{\X}(v_j)\big) _\rho \ \longrightarrow \ \big(\Hc^i(\overline{\X}(c))_\rho\big)^{\oplus 2} \ \longrightarrow \cdots $}

\noindent together with the following properties

\begin{itemize} 

\item  the cuspidal part of  $\Hc^\bullet(\overline{\X}(v_j))$ can be explicitely computed using \hyperref[cusp3]{$\mathrm{(C3)}$}:

\centers{$ \underline{\H}_{\overline{\X}(v_j)} \, = \, (h^8t^4 + h^{10}t^5) \big( \mathrm{E}_8[-\theta]
+\mathrm{E}_8[-\theta^2] +\mathrm{E}_8[\zeta] + \mathrm{E}_8[\zeta^2] + \mathrm{E}_8[\zeta^3] + \mathrm{E}_8[\zeta^4] \big)$}

\item the cuspidal part of a variety associated to a Coxeter element is given by \cite{Lu} (or equivalently can be computed using  \hyperref[cusp3]{$\mathrm{(C3)}$}):

\centers{$ \underline{\H}_{\overline{\X}(c)} \, = \, h^8t^4 \big( \mathrm{E}_8[-\theta]
+\mathrm{E}_8[-\theta^2] +\mathrm{E}_8[\zeta] + \mathrm{E}_8[\zeta^2] + \mathrm{E}_8[\zeta^3] + \mathrm{E}_8[\zeta^4] \big)$}

\item the cohomology of $\X(v_j)$ vanishes in degree $8$.

\end{itemize}

\noindent We obtain, for $j = 1,2,3$:

\centers{$ \underline{\H}_{\X(v_j)} \, = \, (h^9t^4 + h^{10}t^5) \big( \mathrm{E}_8[-\theta]
+\mathrm{E}_8[-\theta^2] +\mathrm{E}_8[\zeta] + \mathrm{E}_8[\zeta^2] + \mathrm{E}_8[\zeta^3] + \mathrm{E}_8[\zeta^4] \big).$}

\noindent Using \cite[Proposition 3.2.10]{DMR}, one can check that it is also the cuspidal part of the cohomology of $\X(\mathbf{v}_4)$. 

\sk

We claim that we can derive the cohomology of $\Z$: for any cuspidal character $\rho$, we have an exact sequence

 \centers{ $\cdots \longrightarrow \displaystyle \bigoplus_{j=1}^4  \Hc^i\big(\X(\mathbf{v}_j)\big)_\rho \ \longrightarrow \Hc^i\big(\Z \big) _\rho \ \longrightarrow \ \big(\Hc^i(\X(c))_\rho\big)^{\oplus 4} \ \longrightarrow \cdots $}

\noindent Furthermore,  the cohomology of $\X(w)$ vanishes in degree $9$ and the cuspidal part of $\Hc^\bullet(\overline{\X}(w))$ is concentrated in degree $10$, given by

\centers{$ h^{10} t^5 \big(\mathrm{E}_8[\mathrm{i}] + \mathrm{E}_8[-\mathrm{i}]  +  3(\mathrm{E}_8[-\theta] +\mathrm{E}_8[-\theta^2]) + 4(\mathrm{E}_8[\zeta] + \mathrm{E}_8[\zeta^2] + \mathrm{E}_8[\zeta^3] + \mathrm{E}_8[\zeta^4] )\big)$} 

\noindent Consequently, the cuspidal part of $\Hc^8(\Z)$ is zero by \ref{5eq10} and we obtain

\centers{$ \underline{\H}_{\Z} \, = \, 4 h^{10} t^5 \big( \mathrm{E}_8[-\theta]
+\mathrm{E}_8[-\theta^2] +\mathrm{E}_8[\zeta] + \mathrm{E}_8[\zeta^2] + \mathrm{E}_8[\zeta^3] + \mathrm{E}_8[\zeta^4] \big).$}

\noindent In particular, we can unpack the exact sequence  \ref{5eq10} according to the different cuspidal characters as follows

\centers{$ 0 \longrightarrow \Hc^{10}\big(\X(w)\big)_{\mathrm{E}_8[\pm\mathrm{i}]} \ \longrightarrow t^5 \mathrm{E}_8[\pm\mathrm{i}] \ \longrightarrow 0 $}


\centers{$ \hskip-1.7mm 0 \longrightarrow \Hc^{10}\big(\X(w)\big)_{\mathrm{E}_8[- \theta^i]} \longrightarrow 3 t^5 \mathrm{E}_8[-\theta^i]  \longrightarrow  4 t^5 \mathrm{E}_8[-\theta^i]    \longrightarrow  \Hc^{11}\big(\X(w)\big)_{\mathrm{E}_8[- \theta^i]} \longrightarrow 0$}

\centers{$ 0 \longrightarrow \ \Hc^{10}\big(\X(w)\big)_{\mathrm{E}_8[\zeta^j]} \ \longrightarrow \ 4 t^5 \mathrm{E}_8[\zeta^j] \  \longrightarrow\  4 t^5 \mathrm{E}_8[\zeta^j]     \ \longrightarrow  \ \Hc^{11}\big(\X(w)\big)_{\mathrm{E}_8[\zeta^j]}\ \longrightarrow 0$}

\noindent To conclude, we observe that the unipotent characters $
\mathrm{E}_8[-\theta^i]$ and  $ \mathrm{E}_8[\zeta^j]$ already appear in the Coxeter variety, and for that reason they cannot be constituents of $\Hc^{10}(\X(w))$ with an eigenvalue of absolute value $q^5$ (see  \hyperref[cusp4]{$\mathrm{(C4)}$}).

\begin{prop}\label{5prop13}Let $w =  t_8 t_7 t_6 t_5 t_4 t_5 t_2 t_4  t_3 t_1 $. The cuspidal part of the  cohomology of $\X(w)$ is given by

\centers{$ h^{10} t^5 \big(\mathrm{E}_8[\mathrm{i}] + \mathrm{E}_8[-\mathrm{i}] \big) +  h^{11} t^{5} \big( \mathrm{E}_8[-\theta] +  \mathrm{E}_8[-\theta^2]\big).$}

\end{prop}

\mk

\noindent \textbf{\thesubsection.3. Cohomology of $\X(w)$.} We summarize the results obtained in this section. The unipotent characters in the principal $\Phi_{24}$-bloc $b$ are given by

\centers{$\begin{array}{r} b_{\mathrm{uni}} \, = \, \big\{ \mathrm{Id}_G, \mathrm{St}_G, \phi_{35,74}, \phi_{160,55}, \phi_{350,38}, \phi_{448,25}, \phi_{350,14}, \phi_{160,16}, \phi_{35,2}, D_{4,\phi_{2,16}''}, \hphantom{\}} \\[5pt] 
D_{4,\phi_{8,9}''} , D_{4,\phi_{12,4}} ,  D_{4,\phi_{8,3}'} ,D_{4,\phi_{2,4}'}, E_6[\theta]_{\phi_{1,3}'},  E_6[\theta]_{\phi_{2,2}}, E_6[\theta]_{\phi_{1,3}''}, E_6[\theta^2]_{\phi_{1,3}'},\hphantom{\}}  \\[5pt]
  E_6[\theta^2]_{\phi_{2,2}}, E_6[\theta^2]_{\phi_{1,3}''},  E_8[\mathrm{i}], E_8[-\mathrm{i}] ,E_8[-\theta], E_8[-\theta^2]  \big\} 
\end{array}$}

\noindent By comparing the restriction to $E_7$ of these characters and Proposition \ref{5prop11}, \ref{5prop12} and \ref{5prop13} we obtain a good approximation of the cohomology of $\X(w)$. 
 
\begin{thm}\label{5thm4}Let $(\G,F)$ be a split group of type $E_8$ and $w$ be a good $24$-regular element of $W$. The contribution to the cohomology of the Deligne-Lusztig $\X(w)$ of the principal series, the $D_4$-series, the $E_6$-series and the cuspidal characters coincides with the contribution of the principal $\Phi_{24}$-block, and it is given by

\centers{\begin{tabular}{l} \begin{tabular}{@{{\vrule width 1pt}\,\ }c@{\,\ {\vrule width 1pt}\,\ }c|c|c|c@{\,\ {\vrule width 1pt}}}\hlinewd{1pt} 
$i $& $\vphantom{\Big(}  10 $& $11 $& $12 $& $13$ \\\hlinewd{1pt} 
$b\H^i(\X(w),\qlb) $&$ \vphantom{\Big(^y} \mathrm{St} $& $q^2 \phi_{35,74} $& $q^3 \phi_{160,55} $& $q^4 \phi_{350,38}   $\\[5pt]
&$\color{purple}-q^3 D_{4,\phi_{2,16}''} $&$\color{purple} -q^4D_{4,\phi_{8,9}''} $&$\color{purple} -q^5 D_{4,\phi_{12,4}} 
$& $\color{purple}-q^6 D_{4,\phi_{8,3}'} $   \\[5pt]
& $\color{violet} \theta q^4 E_6[\theta]_{\phi_{1,3}'} $ & $\color{violet}\theta q^5 E_6[\theta]_{\phi_{2,2}} $ &$\color{violet} \theta q^6 E_6[\theta]_{\phi_{1,3}''} $&  \\[5pt]
& $\color{violet}\theta^2 q^4 E_6[\theta^2]_{\phi_{1,3}'}$ & $\color{violet}\theta^2 q^5 E_6[\theta^2]_{\phi_{2,2}}$ & $\color{violet} \theta^2 q^6 E_6[\theta^2]_{\phi_{1,3}''} $&  \\[5pt]
& $\mathrm{i}q^5 E_8[\mathrm{i}] $& & &  \\[5pt]
& $-\mathrm{i}q^5 E_8[-\mathrm{i}]$ & & &  \\[5pt]
& & $-\theta q^{5} E_8[-\theta ] $& &  \\[5pt]
& & $-\theta^2 q^{5} E_8[- \theta^2] $& &  \\[5pt] \hlinewd{1pt}
\end{tabular}\\ 
\\
\begin{tabular}{@{{\vrule width 1pt}\,\ }c@{\,\ {\vrule width 1pt}\,\ }c|c|c|c|c|c|c@{\,\ {\vrule width 1pt}}} \hlinewd{1pt} $i $& $\vphantom{\Big(}   14 $& $15 $& $16 $& $17$ & $18$ & $19 $& $20$\\\hlinewd{1pt} 
$b\H^i(\X(w),\overline{\mathbb{Q}}_\ell) $& $\vphantom{\Big(^y} q^5 \phi_{448,25} $& $q^6 \phi_{350,14} $& $q^7 \phi_{160,7} $& $q^8 \phi_{35,2} $&  & & $q^{10}\mathrm{Id}  $\\[5pt]
&$\color{purple} -q^7 D_{4,\phi_{2,4}'} $ & & & & & &  \\[5pt] \hlinewd{1pt}
\end{tabular} \end{tabular}}

\noindent where the $D_4$-series is given under Assumption \ref{5hyp1} and the  $E_6$-series under Assumption \ref{5hyp2}.

\end{thm}

\section{Conjectures on associated Brauer trees}

Having computed the cohomology of some Deligne-Lusztig varieties for exceptionals groups, we would like to propose conjectures on Brauer trees for the corresponding principal $\Phi_d$-blocks.

\sk

Recall from \cite{BMM} that if $d$ is a regular number, and $w$ is a $d$-regular element,  the irreducible constituent of the virtual character $\mathrm{R}_{T_w}^G(1) = \sum (-1)^i \Hc(\X(w),\qlb)$ are exactly the unipotent characters in the principal $\Phi_d$-block. If moreover $C_W(wF) \simeq N_G(\T_w)/C_G(\T_w) $ is cyclic, then the $\Phi_d$-block is generically of cyclic defect: if $\ell$ divides $\Phi_d(q)$ but does not divide $ |W|$, then any Sylow subgroup of $G$ is cyclic. In that case, the representation theory of the block (i.e. the module category over the block) can be decribed by its Brauer tree. More precisely, in this situation:

\begin{itemize}

\item  any $\ell$-character $\theta$ of $T_w$ is in general position and the associated irreducible character $\chi_\theta = (-1)^{\ell(w)} \mathrm{R}_{T_w}^G(\theta)$ is cuspidal by \cite[Proposition 2.18]{Lu2}. Moreover, using  \cite[Proposition 12.2]{DM} it can  be shown that its restriction to the set of $\ell$-regular elements is independent from $\theta$. Any character of this form is said to be a \emph{exceptional};

\item there are $e = |C_W(w)|$ unipotent characters $\{\chi_0, \ldots,\chi_{e-1}\}$ in the block, which will be refered as the \emph{non-exceptional} characters;

\end{itemize}

\noindent Now if we consider the sum $\chi_{\mathrm{exc}}$ of all distinct unipotent characters, any projective indecomposable $\overline{\mathbb{F}}_\ell G$-module lifts uniquely, up to isomorphism, to a $\overline{\mathbb{Z}}_\ell$-module $P$ whose character is $[P] = \chi + \chi'$ for $\chi,\chi'$ two distinct characters in $\mathscr{V} = \{\chi_{\mathrm{exc}}, \chi_0, \ldots, \chi_{e-1}\}$. We  define the \emph{Brauer tree} $\Gamma$ of the block to be the graph with vertices labelled by $\mathscr{V}$ and egdes $\chi \, $---$ \, \chi'$ whenever there exists a projective indecomposable module with character $\chi + \chi'$. This graph is a tree and its planar embbeding determines the module category over the block up to Morita equivalence.

\sk

When $d = h$ is a the Coxeter number, Hiss, L\"ubeck and Malle have formulated in \cite{HLM} a conjecture relating the cohomology of the Deligne-Lusztig variety  associated to a Coxeter element (together with the action of $F$) and the planar embedded Brauer tree of the principal $\Phi_h$-block. Using the explicit results on the cohomology of Deligne-Lusztig varieties that we have obtained, and the Brauer trees that we already already know from \cite{HL} and \cite{HLM}, we shall propose two conjectural Brauer trees for groups of type $E_7$ and $E_8$.

\subsection{Observations}

Let $(\G,F)$ be a split group of type $F_4$ and $w$ be a good $8$-regular element. When $\ell$ divides $\Phi_8(q)$ and does not divide the order of $W$, we can observe that the classes in $\F_\ell$ of the eigenvalues of $F$ on $b\Hc^\bullet(\X(w),\qlb)$ form the group of $8$th roots of unity, generated by the class of $q$. Therefore to any non-exceptional character $\chi$ one can associate an integer $j_\chi$ such that the class of the corresponding eigenvalue of $F$ coincides with the class of $q^{j_\chi}$. By  \cite{HL}, the Brauer tree of the block, together with the integers $j_\chi$ is given by

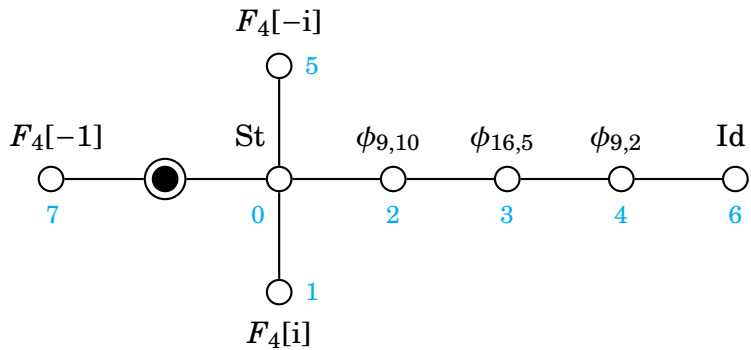
\begin{figure}[h!] 
\centers{\begin{pspicture}(9,4.1)
 \cnode[fillstyle=solid,fillcolor=black](1.5,2){5pt}{A2}
    \cnode(1.5,2){8pt}{A}
  \cnode(3,2){5pt}{B}
  \cnode(4.5,2){5pt}{C}
  \cnode(6,2){5pt}{D}
  \cnode(7.5,2){5pt}{E}
  \cnode(9,2){5pt}{F}
  \cnode(0,2){5pt}{G}
  \cnode(3,0.5){5pt}{I}
  \cnode(3,3.5){5pt}{J}
    \ncline[nodesep=0pt]{G}{A}\naput[npos=-0.1]{$\vphantom{\Big(}F_4[-1]$}
    \nbput[npos=-0.15]{\footnotesize \textcolor{cyan}{$\vphantom{\Big(} 7$}}
    \ncline[nodesep=0pt]{A}{B}\naput[npos=0.8]{$\vphantom{\Big(} \mathrm{St}$}
    \nbput[npos=0.9]{\footnotesize \textcolor{cyan}{$\vphantom{\Big(}0$}}
  \ncline[nodesep=0pt]{B}{C}\naput[npos=1.1]{$\vphantom{\Big(}\phi_{9,10}$}
  \nbput[npos=1.15]{\footnotesize \textcolor{cyan}{$\vphantom{\Big(}2$}}
  \ncline[nodesep=0pt]{C}{D}\naput[npos=1.1]{$\vphantom{\Big(}  \phi_{16,5}$}
  \nbput[npos=1.15]{\footnotesize \textcolor{cyan}{$\vphantom{\Big(}3$}}
  \ncline[nodesep=0pt]{D}{E}\naput[npos=1.1]{$\vphantom{\Big(}\phi_{9,2}$}
  \nbput[npos=1.15]{\footnotesize \textcolor{cyan}{$\vphantom{\Big(}4$}}
  \ncline[nodesep=0pt]{E}{F}\naput[npos=1.1]{$\vphantom{\Big(} \mathrm{Id}$}
  \nbput[npos=1.15]{\footnotesize \textcolor{cyan}{$\vphantom{\Big(}6$}}
  \ncline[nodesep=0pt]{B}{I}\ncput[npos=1.65]{$\vphantom{\Big(} F_4[\mathrm{i}]$}
  \naput[npos=1.15]{\footnotesize \textcolor{cyan}{$\hphantom{i}1$}}
  \ncline[nodesep=0pt]{B}{J}\ncput[npos=1.65]{$\vphantom{\Big(} F_4[-\mathrm{i}]$}
  \nbput[npos=1.15]{\footnotesize \textcolor{cyan}{$\hphantom{i}5$}}
\end{pspicture}}
\caption{Brauer tree of the principal $\Phi_8$-block of $F_4$}
\label{F4}
\end{figure}

Now assume that $(\G,F)$ is a split group of type $E_6$. The Brauer tree of the principal $\Phi_9$-block of $G$ has been determined in \cite{HLM}. It corresponds to the following picture:

\begin{figure}[h!] 
\centers{\begin{pspicture}(10.5,4)

  \cnode[fillstyle=solid,fillcolor=black](0,2){5pt}{A2}
    \cnode(0,2){8pt}{A}
  \cnode(1.5,2){5pt}{B}
  \cnode(3,2){5pt}{C}
  \cnode(4.5,2){5pt}{D}
  \cnode(6,2){5pt}{E}
  \cnode(7.5,2){5pt}{F}
  \cnode(9,2){5pt}{G}
  \cnode(10.5,2){5pt}{H}
  \cnode(1.5,0.5){5pt}{I}
  \cnode(1.5,3.5){5pt}{J}
    \ncline[nodesep=0pt]{A}{B}\naput[npos=0.8]{$\vphantom{\Big(} \mathrm{St}$}
        \nbput[npos=0.9]{\footnotesize \textcolor{cyan}{$\vphantom{\Big(}0$}}
  \ncline[nodesep=0pt]{B}{C}\naput[npos=1.1]{$\vphantom{\Big(}\phi_{20,20}$}
    \nbput[npos=1.15]{\footnotesize \textcolor{cyan}{$\vphantom{\Big(}2$}}
  \ncline[nodesep=0pt]{C}{D}\naput[npos=1.1]{$\vphantom{\Big(}  \phi_{64,13}$}
   \nbput[npos=1.15]{\footnotesize \textcolor{cyan}{$\vphantom{\Big(}3$}}
  \ncline[nodesep=0pt]{D}{E}\naput[npos=1.1]{$\vphantom{\Big(}\phi_{90,8}$}
    \nbput[npos=1.15]{\footnotesize \textcolor{cyan}{$\vphantom{\Big(}4$}}
  \ncline[nodesep=0pt]{E}{F}\naput[npos=1.1]{$\vphantom{\Big(} \phi_{64,4}$}
  \nbput[npos=1.15]{\footnotesize \textcolor{cyan}{$\vphantom{\Big(}5$}}
  \ncline[nodesep=0pt]{F}{G}\naput[npos=1.1]{$\vphantom{\Big(} \phi_{20,2}$}
  \nbput[npos=1.15]{\footnotesize \textcolor{cyan}{$\vphantom{\Big(}6$}}
  \ncline[nodesep=0pt]{G}{H}\naput[npos=1.1]{$\vphantom{\Big(} \mathrm{Id}$}
  \nbput[npos=1.15]{\footnotesize \textcolor{cyan}{$\vphantom{\Big(}8$}}
  \ncline[nodesep=0pt]{B}{I}\ncput[npos=1.65]{$\vphantom{\Big(} E_6[\theta]$}
    \naput[npos=1.15]{\footnotesize \textcolor{cyan}{$\hphantom{i}1$}}
  \ncline[nodesep=0pt]{B}{J}\ncput[npos=1.65]{$\vphantom{\Big(} E_6[\theta^2]$}
    \nbput[npos=1.15]{\footnotesize \textcolor{cyan}{$\hphantom{i}7$}}
\end{pspicture}}
\caption{Brauer tree of the principal $\Phi_9$-block of $E_6$}
\label{E6}
\end{figure}
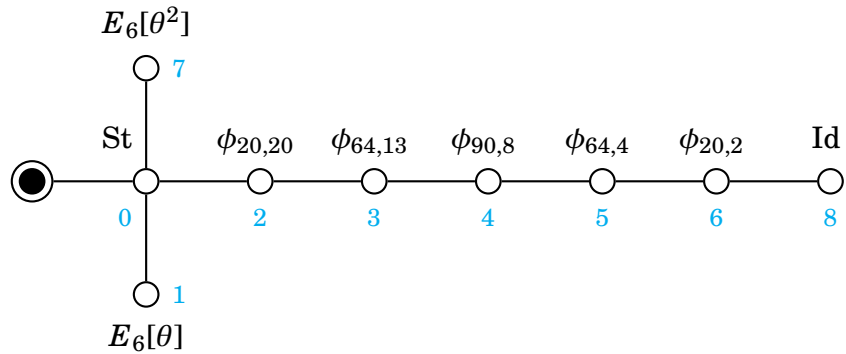
\begin{rmk}\label{5rm20}Unlike the Coxeter case (see \cite{Du3} and \cite{Du4}), the cohomology of the Deligne-Lusztig variety $\widetilde \X(w)$ with coefficients in $\mathbb{Z}_\ell$ is not torsion-free. Indeed, it is impossible to represent the generalized $(q^2)$-eigenspace of $F$ on $\widetilde \X(w)$ with a complex of projective modules $ 0 \longrightarrow P \mathop{\longrightarrow}\limits^f Q\longrightarrow 0$ where the cokernel of $f$ is torsion-free. Note that even the cohomology of the complex constructed by Rickard in \cite[Section 4]{Ri} will also have a non-trivial torsion part (one can show nevertheless that the torsion is always cuspidal). 
\end{rmk}

\subsection{Conjectures}

From the results obtained in Theorem   \ref{5thm3} and \ref{5thm4}, it is not difficult to extrapolate the previous trees to the case of $E_7$ and $E_8$. We conjecture that the Brauer trees of the principal $\Phi_{14}$-block in $E_7$ and the principal $\Phi_{24}$-block in $E_8$ are given by Figure \ref{E7} and \ref{E8}. Note that 

\begin{itemize} 

\item the lines represented by each Harish-Chandra series, as well as the real steam, are known from \cite{Ge3};

\item the simple modules corresponding to edges connecting two different series are necessarily cuspidal.

\end{itemize}

\begin{landscape}

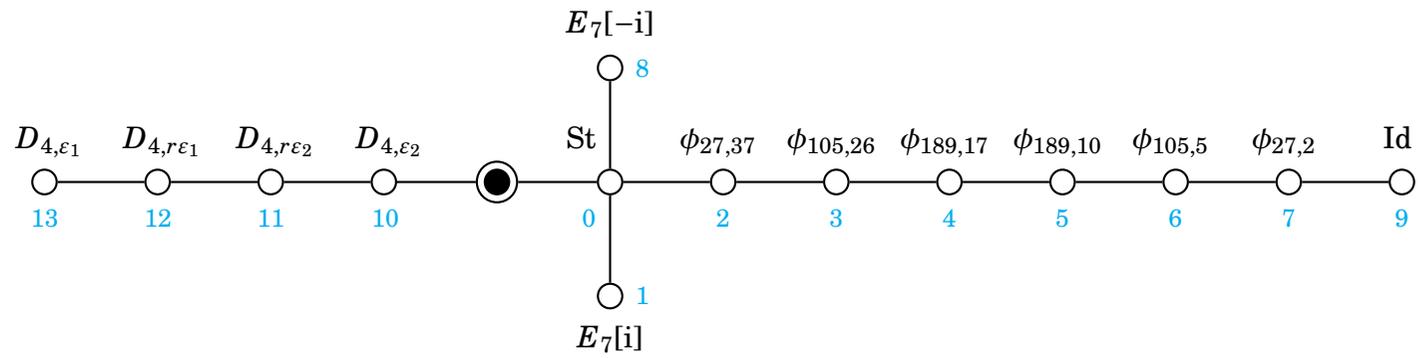
\begin{figure}[h!] \hskip 1cm \begin{pspicture}(18,9)

  \cnode[fillstyle=solid,fillcolor=black](6,2){5pt}{A2}
    \cnode(6,2){8pt}{A}
  \cnode(7.5,2){5pt}{B}
  \cnode(9,2){5pt}{C}
  \cnode(10.5,2){5pt}{D}
  \cnode(12,2){5pt}{E}
  \cnode(13.5,2){5pt}{F}
  \cnode(15,2){5pt}{G}
  \cnode(16.5,2){5pt}{H}
  \cnode(18,2){5pt}{K}  
  \cnode(7.5,0.5){5pt}{I}
  \cnode(7.5,3.5){5pt}{J}
  \cnode(0,2){5pt}{L}
  \cnode(1.5,2){5pt}{M}
  \cnode(3,2){5pt}{N}
  \cnode(4.5,2){5pt}{O}

  \ncline[nodesep=0pt]{A}{B}\naput[npos=0.8]{$\vphantom{\Big(} \mathrm{St}$}
  \nbput[npos=0.9]{\footnotesize \textcolor{cyan}{$\vphantom{\Big(}0$}}
  \ncline[nodesep=0pt]{B}{C}\naput[npos=1.1]{$\vphantom{\Big(}\phi_{27,37}$}
  \nbput[npos=1.15]{\footnotesize \textcolor{cyan}{$\vphantom{\Big(}2$}}
  \ncline[nodesep=0pt]{C}{D}\naput[npos=1.1]{$\vphantom{\Big(}  \phi_{105,26}$}
  \nbput[npos=1.15]{\footnotesize \textcolor{cyan}{$\vphantom{\Big(}3$}}
  \ncline[nodesep=0pt]{D}{E}\naput[npos=1.1]{$\vphantom{\Big(}\phi_{189,17}$}
  \nbput[npos=1.15]{\footnotesize \textcolor{cyan}{$\vphantom{\Big(}4$}}
  \ncline[nodesep=0pt]{E}{F}\naput[npos=1.1]{$\vphantom{\Big(} \phi_{189,10}$}
  \nbput[npos=1.15]{\footnotesize \textcolor{cyan}{$\vphantom{\Big(}5$}}
  \ncline[nodesep=0pt]{F}{G}\naput[npos=1.1]{$\vphantom{\Big(} \phi_{105,5}$}
  \nbput[npos=1.15]{\footnotesize \textcolor{cyan}{$\vphantom{\Big(}6$}}
  \ncline[nodesep=0pt]{G}{H}\naput[npos=1.1]{$\vphantom{\Big(} \phi_{27,2}$}
  \nbput[npos=1.15]{\footnotesize \textcolor{cyan}{$\vphantom{\Big(}7$}}
  \ncline[nodesep=0pt]{B}{I}\ncput[npos=1.65]{$\vphantom{\Big(} E_7[\mathrm{i}] $}
    \naput[npos=1.15]{\footnotesize \textcolor{cyan}{$\hphantom{i}1$}}
  \ncline[nodesep=0pt]{B}{J}\ncput[npos=1.65]{$\vphantom{\Big(} E_7[-\mathrm{i}] $}
    \nbput[npos=1.15]{\footnotesize \textcolor{cyan}{$\hphantom{i}8$}}
  \ncline[nodesep=0pt]{H}{K}\naput[npos=1.1]{$\vphantom{\Big(} \mathrm{Id} $}
  \nbput[npos=1.15]{\footnotesize \textcolor{cyan}{$\vphantom{\Big(}9$}}
  \ncline[nodesep=0pt]{L}{M}\naput[npos=-0.1]{$\vphantom{\Big(} D_{4,\varepsilon_1}$}
  \nbput[npos=-0.15]{\footnotesize \textcolor{cyan}{$\vphantom{\Big(}13$}}
  \ncline[nodesep=0pt]{M}{N}\naput[npos=-0.1]{$\vphantom{\Big(}D_{4,r\varepsilon_1} $}
  \nbput[npos=-0.15]{\footnotesize \textcolor{cyan}{$\vphantom{\Big(}12$}}
  \ncline[nodesep=0pt]{N}{O}\naput[npos=-0.1]{$\vphantom{\Big(}D_{4,r\varepsilon_2} $}
  \nbput[npos=-0.15]{\footnotesize \textcolor{cyan}{$\vphantom{\Big(}11$}}
  \ncline[nodesep=0pt]{O}{A}\naput[npos=-0.1]{$\vphantom{\Big(} D_{4,\varepsilon_2} $}
  \nbput[npos=-0.15]{\footnotesize \textcolor{cyan}{$\vphantom{\Big(}10$}}

\end{pspicture}

\vskip 1cm

\caption{Brauer tree of the principal $\Phi_{14}$-block of $E_7$}
\label{E7}

\end{figure}

\end{landscape}

\begin{landscape}

\begin{figure}[h!] \hskip 1cm\begin{pspicture}(19.5,8)

  \cnode[fillstyle=solid,fillcolor=black](6,2){5pt}{A2}
  \cnode(6,2){8pt}{A}
  \cnode(7.5,2){5pt}{B}
  \cnode(9,2){5pt}{C}
  \cnode(10.5,2){5pt}{D}
  \cnode(12,2){5pt}{E}
  \cnode(13.5,2){5pt}{F}
  \cnode(15,2){5pt}{G}
  \cnode(16.5,2){5pt}{H}
  \cnode(18,2){5pt}{K}  
  \cnode(8.3,0.9){5pt}{I}
  \cnode(8.3,3.1){5pt}{J}
  \cnode(0,2){5pt}{L}
  \cnode(1.5,2){5pt}{M}
  \cnode(3,2){5pt}{N}
  \cnode(4.5,2){5pt}{O}
  \cnode(-1.5,2){5pt}{P}
  \cnode(19.5,2){5pt}{Q}
  \cnode(6.9,3.1){5pt}{R}
  \cnode(6.9,0.9){5pt}{S}
  \cnode(4.9,3.1){5pt}{T}
  \cnode(4.9,0.9){5pt}{U}
  \cnode(3.8,-0.2){5pt}{V}
  \cnode(3.8,4.2){5pt}{W}
  \cnode(2.7,-1.3){5pt}{X}
  \cnode(2.7,5.3){5pt}{Y}

   \ncline[nodesep=0pt]{A}{B}\naput[npos=1]{$\vphantom{\Big(} \mathrm{St}$}
  \nbput[npos=1.15]{\footnotesize \textcolor{cyan}{$\vphantom{\Big(}0$}}
  \ncline[nodesep=0pt]{B}{C}\naput[npos=1.1]{$\vphantom{\Big(}\phi_{35,74}$}
  \nbput[npos=1.35]{\footnotesize \textcolor{cyan}{$\vphantom{\Big(}2$}}
  \ncline[nodesep=0pt]{C}{D}\naput[npos=1.1]{$\vphantom{\Big(}  \phi_{160,55}$}
  \nbput[npos=1.15]{\footnotesize \textcolor{cyan}{$\vphantom{\Big(}3$}}
  \ncline[nodesep=0pt]{D}{E}\naput[npos=1.1]{$\vphantom{\Big(}\phi_{350,38}$}
  \nbput[npos=1.15]{\footnotesize \textcolor{cyan}{$\vphantom{\Big(}4$}}
  \ncline[nodesep=0pt]{E}{F}\naput[npos=1.1]{$\vphantom{\Big(} \phi_{448,25}$}
  \nbput[npos=1.15]{\footnotesize \textcolor{cyan}{$\vphantom{\Big(}5$}}
  \ncline[nodesep=0pt]{F}{G}\naput[npos=1.1]{$\vphantom{\Big(} \phi_{350,14}$}
  \nbput[npos=1.15]{\footnotesize \textcolor{cyan}{$\vphantom{\Big(}6$}}
  \ncline[nodesep=0pt]{G}{H}\naput[npos=1.1]{$\vphantom{\Big(} \phi_{160,7}$}
  \nbput[npos=1.15]{\footnotesize \textcolor{cyan}{$\vphantom{\Big(}7$}}
  \ncline[nodesep=0pt]{B}{I}\ncput[npos=1.70]{$\vphantom{\Big(} E_8[-\theta] $}
  \naput[npos=1.45]{\footnotesize \textcolor{cyan}{$\vphantom{\Big(}1$}}
  \ncline[nodesep=0pt]{B}{J}\ncput[npos=1.70]{$\vphantom{\Big(} E_8[-\theta^2] $}
  \nbput[npos=1.45]{\footnotesize \textcolor{cyan}{$\vphantom{\Big(}9$}}
  \ncline[nodesep=0pt]{H}{K}\naput[npos=1.1]{$\vphantom{\Big(} \phi_{35,2} $}
  \nbput[npos=1.15]{\footnotesize \textcolor{cyan}{$\vphantom{\Big(}8$}}
  \ncline[nodesep=0pt]{L}{M}\naput[npos=-0.1]{$\vphantom{\Big(} D_{4,\phi_{8,3}'} $}  
  \nbput[npos=-0.15]{\footnotesize \textcolor{cyan}{$\vphantom{\Big(}18$}}
  \ncline[nodesep=0pt]{M}{N}\naput[npos=-0.1]{$\vphantom{\Big(} D_{4,\phi_{12,4}} $}
  \nbput[npos=-0.15]{\footnotesize \textcolor{cyan}{$\vphantom{\Big(}17$}}
  \ncline[nodesep=0pt]{N}{O}\naput[npos=-0.1]{$\vphantom{\Big(}D_{4,\phi_{8,9}''} $}
  \nbput[npos=-0.15]{\footnotesize \textcolor{cyan}{$\vphantom{\Big(}16$}}
  \ncline[nodesep=0pt]{O}{A}\naput[npos=-0.1]{$\vphantom{\Big(}D_{4,\phi_{2,16}''}  $}
  \nbput[npos=-0.15]{\footnotesize \textcolor{cyan}{$\vphantom{\Big(}15$}}
  \ncline[nodesep=0pt]{P}{L}\naput[npos=-0.1]{$\vphantom{\Big(}D_{4,\phi_{2,4}'}  $}
  \nbput[npos=-0.15]{\footnotesize \textcolor{cyan}{$\vphantom{\Big(}19$}}
  \ncline[nodesep=0pt]{K}{Q}\naput[npos=1.1]{$\vphantom{\Big(} \mathrm{Id} $}
  \nbput[npos=1.15]{\footnotesize \textcolor{cyan}{$\vphantom{\Big(}10$}}
  \ncline[nodesep=0pt]{R}{A}\ncput[npos=-0.75]{$\phantom{\Big(} E_8[\mathrm{i}] $}
  \ncline[nodesep=0pt]{S}{A}\ncput[npos=-0.75]{$\phantom{\Big(} E_8[-\mathrm{i}] $}
  \ncline[nodesep=0pt]{T}{A}\naput[npos=-0.1]{$ E_6[\theta]_{\phi_{1,3}'} $}
  \nbput[npos=-0.65]{\footnotesize \textcolor{cyan}{$\vphantom{\Big(}12$}}
  \naput[npos=0.7]{\footnotesize \textcolor{cyan}{$\phantom{\Big(a}11$}}
  \ncline[nodesep=0pt]{U}{A}\nbput[npos=-0.1]{$ \hskip-0.4mm E_6[\theta^2]_{\phi_{1,3}'} $}
  \naput[npos=-0.65]{\footnotesize \textcolor{cyan}{$\vphantom{\Big(}20$}}
  \nbput[npos=0.7]{\footnotesize \textcolor{cyan}{$\phantom{\Big(a}23$}}
  \ncline[nodesep=0pt]{V}{U}\nbput[npos=-0.1]{$E_6[\theta^2]_{\phi_{2,2}} $}
  \naput[npos=-0.45]{\footnotesize \textcolor{cyan}{$\vphantom{\Big(}21$}}
  \ncline[nodesep=0pt]{W}{T}\naput[npos=-0.1]{$E_6[\theta]_{\phi_{2,2}}$}
  \nbput[npos=-0.45]{\footnotesize \textcolor{cyan}{$\vphantom{\Big(}13$}}
  \ncline[nodesep=0pt]{X}{V}\nbput[npos=-0.1]{$ E_6[\theta^2]_{\phi_{1,3}''} $}
  \naput[npos=-0.25]{\footnotesize \textcolor{cyan}{$\vphantom{\Big(}22$}}
  \ncline[nodesep=0pt]{Y}{W}\naput[npos=-0.1]{$E_6[\theta]_{\phi_{1,3}''} $}
  \nbput[npos=-0.25]{\footnotesize \textcolor{cyan}{$\vphantom{\Big(}14$}}

\end{pspicture}

\vskip 2.5cm

\caption{Brauer tree of the principal $\Phi_{24}$-block of $E_8$}
\label{E8}

\end{figure}

\end{landscape}

\begin{center} \subsection*{Acknowledgements}\end{center}

I am indebted to Jean Michel for introducing me to the methods of \cite{DMR} and \cite{DM2} which I have used throughout this paper. Most of this work was carried out while I was visiting him and François Digne. I would like to thank them for many valuable comments and suggestions. 

\bk

\bibliographystyle{abbrv}
\bibliography{cohotherregelt}

\end{document}